\documentclass[12pt,twoside]{article}
\usepackage{amssymb,pb-diagram,CD,papersaver6k}
\newcommand{\PSbox}[3]{\mbox{\rule{0in}{#3}\includegraphics{#1}\hspace{#2}}} 
\newtheorem{Theorem}{Theorem}[section]

\newtheorem{Lemma}[Theorem]{Lemma}
\newtheorem{Proposition}[Theorem]{Proposition}

\newtheorem{Corollary}[Theorem]{Corollary}
\newtheorem{Problem}[Theorem]{Problem}
\newtheorem{Conjecture}[Theorem]{Conjecture}
\newtheorem{Claim}[Theorem]{Claim}

\newcommand{\lra}{\longrightarrow}
\newcommand{\supp}{{\mathrm s}{\mathrm u}{\mathrm p}{\mathrm p}}
\newcommand{\iint}{{\mathrm i}{\mathrm n}{\mathrm t}}
\def\sqr#1#2{{\vcenter{\vbox{\hrule  height.#2pt
        \hbox{\vrule width.#2pt height#1pt \kern#1pt \vrule width.#2pt}
        \hrule height.#2pt}}}}

\let\phi=\varphi
\def\Box{\sqr66}   
\let\epsilon=\varepsilon
\def\){ \right) }
\def\({ \left( }
\def\[{ \left[ }
\def\]{ \right] }
\def\<{ \langle }
\def\>{ \rangle }
\let\ljunk=\{
\let\rjunk=\}
\def\{{\left\ljunk}
\def\}{\right\rjunk}
\def\p{\partial}

\def\conf{{\mathrm c}{\mathrm o}{\mathrm n}{\mathrm f}}

\def\Riem{{\cal R}{\mathrm i}{\mathrm e}{\mathrm m}}
\def\dist{{\mathrm d}{\mathrm i}{\mathrm s}{\mathrm t}}

\newcommand{\Diff}{{\mathrm D}{\mathrm i}{\mathrm f}{\mathrm f}}

\def\Vol{{\mathrm V}{\mathrm o}{\mathrm l}}
\def\Ric{\mbox{Ric}}
\def\Tr{{\mathrm T}{\mathrm r}}

\newcommand{\R}{{\mathbf R}}

\begin{document} 

\title{Relative Yamabe Invariant}
\author{Kazuo Akutagawa\thanks{\ \ Partialy supported by the \
Grants-in-Aid \ for \ Scientific Research, \ The Ministry of Education,
Science, Sports and Culture, Japan, No. 09640102
\protect\\
email: smkacta\atsign ipc.shizuoka.ac.jp
}
\and
Boris Botvinnik\thanks{
\ \ Partially supported by SFB 478, M\"unster
\protect\\
e-mail: botvinn\atsign math.uoregon.edu
}}

\maketitle
\begin{abstract}
We define a relative Yamabe invariant of a smooth manifold with given
conformal class on its boundary. In the case of empty boundary the
invariant coincides with the classic Yamabe invariant. We develop
approximation technique which leads to gluing theorems of two
manifolds along their boundaries for the relative Yamabe invariant. We
show that there are many examples of manifolds with both positive and
non-positive relative Yamabe invariants. 
\end{abstract}
\vspace{-1.5em}
\section{Introduction}\label{int}
{\bf \ref{int}.1. General setting.} Let $W$ be a compact smooth
manifold with boundary, $\p W = M\neq \emptyset$, and $n=\dim W \geq
3$. Let 
$\Riem(W)$ be the space of all Riemannian metrics on $W$.  For a
metric $\bar{g}\in \Riem(W)$ we denote $H_{\bar{g}}$ the mean
curvature along the boundary $\p W = M$, and $g=\bar{g}|_M$. We also
denote $[\bar{g}]$ and $[g]$ the corresponding conformal classes, and
${\cal C}(M)$ and ${\cal C}(W)$ the space of conformal classes on $M$
and $W$ respectively.  Let $\bar{C}$ and $C$ be conformal classes of
metrics on $W$ and $M$ respectively.  We say that $C$ is the {\it
boundary of $\bar{C}$} or $\bar{C}$ is a {\it coboundary of $C$} if
$\bar{C}|_{M} = C$.  We use notation $\p \bar{C}=C$ in this case. Let
${\cal C}(W,M)$ be the space of pairs $(\bar{C},C)$, such that $\p
\bar{C}=C$. Denote $\bar{C}^0 = \{\bar{g}\in \bar{C} \ | \ H_{\bar{g}}=0 \ \}$.
We call  $\bar{C}^0\subset \bar{C}$ the {\sl normalized conformal class}.
Let ${\cal C}^0(W,M)$ be the space
of pairs $(\bar{C}^0,C)$, so that $\bar{C}^0\subset \bar{C}$, and
$(\bar{C},C)\in {\cal C}(W,M)$. It is easy to observe (see
\cite[formula (1.4)]{Escobar1}) that for any conformal class
$\bar{C}\in {\cal C}(W)$ the subclass $\bar{C}^0$ is not empty. Thus
there is a natural bijection between the spaces ${\cal C}^0(W,M)$ and
${\cal C}(W,M)$.  Let $\bar{g}\in \bar{C}^0$ be a metric. Then
$\bar{C}^0$ may be described as follows:
$$
\bar{C}^0 = \{ u^{{4\over n-2}} \bar{g} \ | \ u\in C^{\infty}_+(W) \
\mbox{such that} \ \ \p_{\nu}u = 0 \ \ \mbox{along} \ M \ \}.
$$
Here $\nu$ is a normal unit (inward) vector field along the boundary, and 
$C^{\infty}_+(W)$ is the space of positive smooth functions on $W$.
\vspace{2mm}

\noindent
{\bf \ref{int}.2. The Einstein-Hilbert functional.} Let $C\in
{\cal C}(M)$ be given. We define the following subspaces of metrics:
$$
\begin{array}{lcl}
\Riem_C(W,M) & = & \{ \bar{g}\in \Riem(W) \ | \  \p [\bar{g}]=C \},
\\
\\
\Riem_C^0(W,M) & = & \{ \bar{g}\in \Riem_C(W) \ | \   H_{\bar{g}} = 0 \}.
\end{array}
$$
We consider the normalized Einstein-Hilbert functional $ I :
\Riem_C^0(W,M) \rightarrow \R $ given by
$$
I(\bar{g})= {\int_W R_{\bar{g}} d\sigma_{\bar{g}}\over
\Vol_{\bar{g}}(W)^{{n-2\over n}}},
$$
where $R_{\bar{g}}$ and is the scalar curvature, and
$d\sigma_{\bar{g}}$ is the volume element.  As in the case of closed
manifolds, we have the following result.
\begin{Theorem}\label{Th1}
Critical points of the functional $I$ on the space $\Riem_C^0(W,M)$
coincide with the set of Einstein metrics $\bar{g}$ on $W$ with $\p
[\bar{g}]=C$, and $H_{\bar{g}} = 0$.
\end{Theorem}
We postpone the proof of Theorem \ref{Th1} to Section \ref{Proof-EH}.
\vspace{2mm}

\noindent
{\bf \ref{int}.3. Relative Yamabe invariants.} Similarly to the case
of closed manifolds, the functional $I$ is not bounded. Precisely, it
is easy to prove that for any manifold $W$, $\p W= M$, $\dim W\geq 3$,
and any conformal class $C\in {\cal C}(M)$
$$
{\inf}_{\bar{g}\in \Riem_C^0(W,M)} I(\bar{g}) = -\infty, \ \ \mbox{and} \ \ 
{\sup}_{\bar{g}\in \Riem_C^0(W,M)} I(\bar{g}) = \infty.
$$
Let $(\bar{C},C)\in {\cal C}(W,M)$.  We define the {\it relative
Yamabe constant of} $(\bar{C},C)$ as 
$$
Y_{\bar{C}}(W,M;C) = \inf_{\bar{g}\in \bar{C}^0} I(\bar{g}).
$$
{\bf Remark.} We notice that the Yamabe constant $Y_{\bar{C}}(W,M;C)$
coincides with the constant $Q(W)$ (up to a universal positive factor
depending only on the dimension of $W$) defined by J. Escobar
\cite{Escobar1} for each pair of conformal classes $(\bar{C},C)\in
{\cal C}(W,M)$.
\vspace{2mm}

\noindent
Clearly the Yamabe constant $Y_{\bar{C}}(W,M;C)$ must be related to
the Yamabe problem on a manifold with boundary, which was solved by
P. Cherrier \cite{Cherrier} and J. Escobar \cite{Escobar1} under some
restrictions. Indeed, P. Cherrier proved the existence of a minimizer
for the Yamabe functional $I|_{\bar{C}^0}$ provided
\begin{equation}\label{eq1-a-1}
Y_{\bar{C}}(W,M;C) < Y_{[\bar{g}_0]}(S^n_+,S^{n-1};[g_0]).
\end{equation}
Here $S^n_+$ is a round hemisphere with standard metric $\bar{g}_0$,
and $S^{n-1}\subset S^n_+$ the equator with
$g_0=\bar{g}_0|_{S^{n-1}}$. More generally, J. Escobar \cite{Escobar1}
solved the Yamabe problem under restrictions we list below. We
emphasize that the Escobar's result includes the case when the
inequality (\ref{eq1-a-1}) is satisfied. Here is the list of
restrictions given in \cite{Escobar1}:
\begin{equation}\label{eq1-a}
\begin{array}{lll}
\mbox{(a)}& n\geq 6 & \\ \mbox{(b)}& M =\p W \ \mbox{is umbilic in
$W$} & \mbox{\hspace*{70mm}} \\ \mbox{(c)}& \mbox{the Weyl tensor
$W_C=0$ on $M$} & \\ \mbox{(d)}& \mbox{the Weyl tensor $W_{\bar{C}}
\not\equiv 0$ on $W$}. &
\end{array}
\end{equation}
Notice that the conditions (\ref{eq1-a}) are conformally invariant.
We denote 
$$
{\cal C}^{Esc}(W,M) = \{ (\bar{C},C)\in {\cal C}(W,M) \ \left|  
\begin{array}{c}
\mbox{at least one of the conditions} \\ \mbox{(a)-(d) from
(\ref{eq1-a}) is not satisfied}
\end{array} \right.
\}
$$
{\bf Remark.} It is easy to see that ${\cal C}^{Esc}(W,M)\subset {\cal
C}(W,M)$ is open dense. 
\vspace{1mm}

\noindent
We state the Escobar's result using the
terms introduced above.
\begin{Theorem}\label{esc}
{\rm (Escobar, \cite[Theorem 6.1]{Escobar1})} Let $(W,M)$ be a
compact manifold with boundary, and $(\bar{C},C)\in {\cal
C}^{Esc}(W,M)$. Then there exists a metric $\check{g}\in \bar{C}^0$ such
that $Y_{\bar{C}}(W,M;C)= I(\check{g})$. Such metric $\check{g}$ is
called a {\rm relative Yamabe metric.}
\end{Theorem}
{\bf Remark.}  A relative Yamabe metric $\check{g}\in \bar{C}^0$, of
course, has constant scalar curvature $R_{\check{g}}=
Y_{\bar{C}}(W,M;C)\cdot \Vol_{\check{g}}(W)^{-{2\over n}}$.
\vspace{1mm}

\noindent
We define the {\it relative Yamabe invariant with respect to a conformal
class $C\in {\cal C}(M)$:}
$$
Y(W,M,C) = \sup_{\bar{C}, \p \bar{C} = C}
Y_{\bar{C}}(W,M;C).
$$
Then we would like to define the following relative Yamabe invariants:
$$
\begin{array}{ccl}
Y(M;C) &=& \displaystyle
\sup_{W, \p W = M} Y(W,M;C),
\\ \vspace{-.5em}
\\
Y(W,M) &=& \displaystyle
\sup_{C\in {\cal C}(M)} Y(W,M;C),
\\ \vspace{-.5em}
\\
{\Bbb Y}(M) &=& \displaystyle
\sup_{W, \p W = M} Y(W;M) .
\end{array}
$$
The invariants $Y(W,M;C)$, $Y(M;C)$ have clear geometric meaning in
terms of positive scalar curvature (abbreviated as psc). We call a
conformal class $C\in {\cal C}(M)$ {\it positive} if the Yamabe
constant $Y_C(M)>0$. Of course, it means that a Yamabe metric $g\in C$
has positive constant scalar curvature. The following statement
follows from the above definitions.
\begin{Claim}\label{Cl5}$\mbox{ \ }$ \\
{\bf (1)} Let $C\in{\cal C}(M)$ be a positive conformal class. Then
\begin{enumerate}
\item[$\bullet$] $Y(W,M;C) > 0$ if and only if any psc-metric
$g\in C$ may be extended conformally to a psc-metric $\bar{g}$ on $W$ with
$H_{\bar{g}}=0$.
\item[$\bullet$] $Y(M;C) > 0$ if and only if for any psc-metric
$g\in C$ there exists a manifold $W$, so that $g$ may be extended
conformally to a
psc-metric $\bar{g}$ on $W$ with $H_{\bar{g}}=0$.  
 \end{enumerate}
{\bf (2)} The invariants $Y(W,M)$, ${\Bbb Y}(M)$ are
diffeomorphism invariants. 
\end{Claim}
We present our main results on the relative Yamabe invariant in the
next section.
\vspace{2mm}

\noindent
{\bf \ref{int}.4. Acknowledgments.} Both authors would like to
acknowledge partial financial support provided by the Department of
Mathematics at the University of Oregon, SFB 478 -- Geometrische
Strukturen in der Mathematik, and the Grants-in-Aid for Scientific
Research, Japan. We would like also to thank Professors Mitsuhiro Itoh
and Thomas Schick for useful discussions and their hospitality while
we visited the University of Tsukuba and during our visit at
M\"unster. The second author would like to thank Michael Joachim,
Wolfgang L\"uck, and the members of the geometry-topology in M\"unster
for very intersting mathematical discussions.

\section{Overview of the results}\label{results}
{\bf \ref{results}.1. Minimal boundary condition and approximation
theorems.} First, one can notice that the minimal boundary condition
$H_{\bar{g}}=0$ is too weak for applications. For instance, to apply a
relative index theory, one needs much stronger condition that a metric
$\bar{g}$ is a product metric near the boundary.
\vspace{2mm}

\noindent
The closest differential-geometric approximation to the product metric
near the boundary is when a boundary is totally geodesic.  In more
detail, let $\bar{g}\in \Riem_C^0(W,M)$, $g = \bar{g}|_M$, and, as
above $A_{\bar{g}}=\(\bar{A}_{ij}\)$ be the second fundamental form of
$M=\p W$ with respect to $\bar{g}$.  The boundary $M=\p W$ is said to
be {\it totally geodesic} if $A_{\bar{g}}$ vanishes identically on
$M$. Clearly any metric from the normalized conformal class
$[\bar{g}]^{0}$ is totally geodesic if $\bar{g}$ is. We call the
conformal class $[\bar{g}]$ of such metric $\bar{g}$ {\sl umbilic}.
We denote ${\cal C}^{um}_C(W,M)\subset {\cal C}_C(W,M) :=\{ \bar{C}\in
{\cal C}(W) \ | \ \bar{C}|_M =C\} $ the subspace of umbilic conformal
classes.
\vspace{2mm}

\noindent
Our first aim is to prove a generalization (Proposition \ref{trick})
of the approximation theorem due to Kobayashi \cite{Kobayashi1}.  We
show that any metric $\bar{g}$ with totally geodesic boundary is
$C^1$-close to a metric $\tilde{g}$ which is conformally equivalent to
a product metric near the boundary. Moreover, we show that the scalar
curvature $R_{\bar{g}}$ is $C^0$-close to $R_{\tilde{g}}$ of
$\tilde{g}$.
\vspace{2mm}

\noindent
Next, we prove the approximation Theorem \ref{main-trick} under the
minimal boundary condition. Theorem \ref{main-trick} gives us a
fundamental tool on the relative Yamabe invariant. In particular, we
prove the following result.
\begin{Theorem}\label{Th2}
For any $\bar{C}\in {\cal C}_C(W,M)$, and any $\epsilon>0$ there
exists a conformal class $\tilde{C}\in {\cal C}^{um}_C(W,M)$, and a
metric $\tilde{g}\in \tilde{C}^{0}$, such that
\begin{equation}
\{
\begin{array}{l}
\bar{C} \ \mbox{and} \ \tilde{C} \ \mbox{are $C^0$-close conformal classes}
\\
|Y_{\bar{C}}(W,M;C) - Y_{\tilde{C}}(W,M;C)| < \epsilon
\\
\tilde{g} \sim g + dr^2 \ \ \mbox{(conformally equivalent near $M$),}
\end{array}
\right.
\end{equation}
where $C=\p\tilde{C}$ and $g = \bar{g}|_M$. More precisely, 
$$
\begin{array}{l}
\tilde{g}=\(1+{r^2\over 2} f(x)\)^{{4\over n-2}} (g+ dr^2) \ \ \ 
\mbox{near $M$, where}
\\
\\
f(x)=-{n-2\over
4(n-1)}(R_{\bar{g}}- R_g + 3|A_{\bar{g}}|_g^2) \ \ \   \mbox{on $M$}.
\end{array}
$$
\end{Theorem}
We define the ``umbilic Yamabe invariant'' $Y^{um}(W,M;C)$ as
$$
Y^{um}(W,M;C)=\sup_{\tilde{C}\in {\cal C}^{um}_C(W,M)}
Y_{\tilde{C}}(W,M;C).
$$
Theorem \ref{Th2} leads to the following conclusion.
\begin{Corollary}
$Y^{um}(W,M;C)= Y(W,M;C)$.
\end{Corollary}
{\bf \ref{results}.2. Gluing Theorem.} We analyze the gluing
procedure for manifolds equipped with conformal structures.  Let
$W_1$, $W_2$ be two compact manifolds, $\dim W_1 = \dim W_2 \geq 3$,
with boundaries
$$
\p W_1 = M_1= M_0 \sqcup M,  \ \ \ \mbox{and}\ \ \ 
\p W_2 = M_2= M_0 \sqcup M^{\prime}
$$
endowed with conformal classes $C_1 = C_0\sqcup C\in {\cal C}(M_1)$,
$C_2 = C_0\sqcup C^{\prime}\in {\cal C}(M_2)$, where $C_0\in {\cal
C}(M_0)$, $C\in {\cal C}(M)$, $C^{\prime}\in {\cal C}(M^{\prime})$.
Let $W= W_1\cup_{M_0}(-W_2)$ be the boundary connected sum 
of $W_1$ and $W_2$ along $M_0$. 
\vspace{2mm}

\noindent
We study the case when the conformal class $C_0\in {\cal C}(M_0)$ is
positive, and the relative Yamabe invariants $Y(W_j,M_j;C_j)$, $j=1,2$
are positive as well. We essentially use the approximation Theorem
\ref{main-trick} to prove the following result (Theorem \ref{ThI}):
\begin{Theorem}\label{ThI-1}
Let $C_0\in {\cal C}(M_0)$ be a positive conformal class, and
$Y(W_j,M_j;C_j)> 0 $ for $j=1,2$. Then $Y(W,\p W; C\sqcup
C^{\prime})>0$.
\end{Theorem}
In particular, this result allows us to construct many examples with
positive relative Yamabe invariants (Theorem \ref{ThII}). 
\vspace{2mm}

\noindent
{\bf \ref{results}.3. Non-positive Yamabe invariant.}  Let $W$ be a
compact smooth $n$-manifold with boundary $\p W= M\neq \emptyset$, and
$M= M_{0} \sqcup M_{1}$.  (Here $M_{1}$ may be empty.) Let $(\bar{C},
C = C_{0} \sqcup C_{1})$ be a pair of conformal classes on $(W, M)$,
where $C_0\in {\cal C}(M_0)$, $C_1\in {\cal C}(M_1)$.  Let $X = W
\cup_{M_{0}}(-W)$ be the double of $W$ along $M_{0}$.  Then, the
boundary of $X$ is $\p X= M_{1} \sqcup (-M_{1})$.  We prove the
following result (Theorem \ref{[Theorem 4.12]}):
\begin{Theorem}\label{[Theorem 4.12]-1}$\mbox{\ }$
\begin{enumerate}
\item[{\bf (1)}] If  $Y_{\bar{C}}(W, M; C)\leq 0$, then

\hspace*{10mm}$ \displaystyle 
2^{{2\over n}}  Y_{\bar{C}}(W, M; C)  \leq  
    Y(X, M_{1} \sqcup (-M_{1}); C_{1}  \sqcup  C_{1}).
$
\item[{\bf (2)}] If  $Y_{\bar{C}}(W, M; C)\leq 0$, then

\hspace*{10mm}$ \displaystyle 
    2^{{2\over n}}  Y(W, M; C)  \leq 
    Y(X, M_{1} \sqcup (-M_{1}); C_{1}  \sqcup C_{1}).  
$
\item[{\bf (3)}]   If  $Y(W, M)\leq 0$, then

\hspace*{10mm}
$ \displaystyle 
2^{{2\over n}} Y(W, M) \leq Y(X, M_{1} \sqcup
(-M_{1})).  $
\end{enumerate}
\end{Theorem}
When the manifolds $\p W= M = M_{0}$, and $M_{1}$ is empty, (in this case,
the boundary of $X = W\cup_{M}(-W)$ is empty)  the following holds
(Corollary \ref{[Corollary 4.12]}):
\begin{Corollary}\label{[Corollary 4.12]-1}$\mbox{\ }$
\begin{enumerate}
\item[{\bf (1)}]   
If  $Y_{\bar{C}}(W, M; C)\leq 0$, \ \ then \ \ 
$
    2^{{2\over n}} Y_{\bar{C}}(W, M; C)  \leq  Y(X).
$ 
\item[{\bf (2)}]   If  $Y(W, M; C)\leq 0$, \ \ \  then \ \
$    
 2^{{2\over n}} Y(W, M; C)  \leq  Y(X). 
$ 
\item[{\bf (3)}]   If  $Y(W, M)\leq 0$, \ \ \ \ \ \ \ then \ \ 
$
   2^{{2\over n}} Y(W, M) \leq Y(X). 
$
\end{enumerate}
\end{Corollary}
We use these results to show that there are many examples of manifolds
with non-positive Yamabe invariant. In particular, we have (Corollary
\ref{enlarg}):
\begin{Corollary}\label{enlarg-1}
Let $N$ be an enlargeable compact closed manifold. Then $$2^{{2\over
n}} Y(N\setminus \iint (D^n), S^{n-1})\leq Y(N\# (-N)).$$ In particular,
$Y(N\setminus \iint (D^n), S^{n-1})\leq 0$.
\end{Corollary}
The rest of the paper is organized as follows. We prove Theorem
\ref{Th1} in Section \ref{Proof-EH}. Then we prove the approximation
theorems in Section \ref{AB-trick}. We give a gluing construction in
Section \ref{diss}. We analyze the Yamabe invariant for a double,
and study the non-positive Yamabe invariant in Section
\ref{non-pos}. In the last Section \ref{open} we define and study the
moduli space of positive conformal classes and introduce conformal
concordance and conformal cobordism.
\section{Proof of Theorem \ref{Th1}}\label{Proof-EH}  Let $\bar{g}\in
\Riem_C(W,M)$ be a metric, and $\{ \bar{g}(t)\}$ be a variation of
$\bar{g}$ in the space $\Riem_C(W,M)$, i.e. $\bar{g}(0)=\bar{g}$. Thus
we consider first a general variation, i.e. $\{ \bar{g}(t)\}$ is not
necessarily in the subspace $\Riem_C^0(W,M)\subset \Riem_C(W,M)$. Now we
need the following notations. Let $ h= \left.{d\over dt}\right|_{t=0}
\bar{g}(t) $ be a variational vector, and
$$
g(t)= \bar{g}(t)|_{M}, \ \ \ \mbox{where} \ \ \ g(0)=g.
$$
{\bf Remark. } We observe that the condition $g(t)\in C$ implies that
$h_{ij}= fg_{ij}$ on $M$, where $f\in C^{\infty}(M)$.
\vspace{2mm}

\noindent
Let $r=r(t)$ be the distance function to the boundary $M$ in $W$ with
respect to the metric $\bar{g}(t)$. Let $\nu= {\p \over \p r}$ be a
unit normal (inward) vector field along the boundary $\p W = M$.
\vspace{2mm}

\noindent
Let $p\in M$, and $\{r,x^1,\ldots,x^{n-1}\}$ be a Fermi coordinate
system near $p$. We use indices $\alpha,\beta=0,1,\ldots, n-1$, where
$0$ correponds to the normal direction, and $i,j,k=1,\ldots,n-1$ are
indices corresponding to the tangent directions (only on the boundary $\p
W = M$). We denote $(\cdot)^{\prime}= \left.{d\over
dt}(\cdot)\right|_{t=0}$, the variational derivative evaluated at
$t=0$.  In order to prove Theorem \ref{Th1}, it is enough to prove the
following formula.
\begin{Claim}\label{Cl1} Let $\{ \bar{g}(t)\}$ be a variation as above.
Then
$$
\(\int_{W} R_{\bar{g}(t)} d\sigma_{\bar{g}(t)}\)^{\prime}
= -\int_{W} \< Ric_{\bar{g}} -{1\over 2} R_{\bar{g}}\bar{g},
 h\>_{\bar{g}} d\sigma_{\bar{g}} -
\int_{M} \( 2 H^{\prime}_{\bar{g}} + f H_{\bar{g}}\) d\sigma_g.
$$
\end{Claim}
{\bf Proof.} We denote $\bar{\nabla}$ and $\nabla$ corresponding
Levi-Civita connections with respect to the metrics $\bar{g}$ and $g$.
Standard calculation gives:
\begin{equation}\label{eq2}
\begin{array}{lcl}
(R_{\bar{g}(t)})^{\prime} & =&
\displaystyle
-\bar{\nabla}^{\alpha}\bar{\nabla}_{\alpha}(\Tr_{\bar{g}} h) + 
\bar{\nabla}^{\alpha}\bar{\nabla}^{\beta} h_{\alpha\beta} - 
\<Ric_{\bar{g}}, h\>_{\bar{g}},
\\
\\
(d\sigma_{\bar{g}(t)})^{\prime} &=&
\<{1\over 2} R_{\bar{g}}\bar{g}, h\>_{\bar{g}}d\sigma_{\bar{g}}.
\end{array}
\end{equation}
The formula (\ref{eq2}) together with Gauss divergence formula gives
$$
\begin{array}{lcl}
\displaystyle
\(\int_{W} R_{\bar{g}(t)} d\sigma_{\bar{g}(t)}\)^{\prime}
& = & \displaystyle
-\int_{W} \< Ric_{\bar{g}} -{1\over 2} R_{\bar{g}}\bar{g},
 h\>_{\bar{g}} d\sigma_{\bar{g}} 
\\
\\
& & \displaystyle
+ \int_{M} \<\bar{\nabla}(\Tr_{\bar{g}}h), \nu\>_{\bar{g}}  d\sigma_g -
\int_{M}\sum_{\alpha=0}^{n-1}
(\bar{\nabla}_{e_{\alpha}}h)(e_{\alpha}, \nu)d\sigma_g.
\end{array}
$$
Here $\{e_{\alpha}\}=\{\nu,e_1,\ldots,e_{n-1}\}$ is a local orthonormal
field. We denote
$$
B_{I} 
=\<\bar{\nabla}(\Tr_{\bar{g}} h ), \nu\>_{\bar{g}} ,
\ \ \ \ \ 
B_{II} 
= - \sum_{\alpha=0}^{n-1}
(\bar{\nabla}_{e_{\alpha}}h)(e_{\alpha}, \nu).
$$
Let $p\in M$ be an arbitrary point of the boundary. As before, let
$\nu$ be a unit vector field normal (inward) to the boundary, such
that $\bar{\nabla}_{\nu}\nu=0$, and $\{e_i\}$ be an orthonormal frame
near $p$ in $W$ such that $\nabla_{e_i }e_j =0$ at $p$ and $t=0$. We
emphasize that, in general, $\bar{\nabla}_{e_i}e_j$ does not vanish at
$p$.  We have the second fundamental form of $M$:
$$
A_{ij} = A(e_i,e_j)= \bar{g}(\bar{\nabla}_{e_i}e_j,\nu) = - 
\bar{g}(\nabla_{e_i}\nu,e_j).
$$
Then we have $H=g^{ij}A_{ij}$ the mean curvature of the boundary $M$.
We have:
$$
H= -\sum_{i=0}^{n-1}\bar{g}(\bar{\nabla}_{e_i}\nu,e_i) = \sum_{\alpha
=0}^{n-1}\bar{g}(\bar{\nabla}_{e_{\alpha}}\nu,e_{\alpha}) =
-\nabla_{\alpha}\nu^{\alpha} = - \( \p_{\alpha}\nu^{\alpha} +
\bar{\Gamma}^{\alpha}_{\alpha\beta}\nu^{\beta}\)
$$
Here $\{x^{\alpha}\}= \{r,x^1,\ldots,x^{n-1}\}$ are a Fermi coordinate
system near $p$ in $W$, and $\displaystyle \p_{\alpha}= {\p\over
\p x^{\alpha}}$ (and $\p_{\alpha}=e_{\alpha}$ at $p$). We have:
$$
\begin{array}{lcl}
H^{\prime} &=& \displaystyle 
- \[ 
\p_{\alpha} (\nu^{\prime})^{\alpha} + \mbox{\begin{picture}(0,15)\end{picture}}
\bar{\Gamma}^{\alpha}_{\alpha\beta}(\nu^{\prime})^{\beta}
+(\bar{\Gamma}^{\prime})^{\alpha}_{\alpha\beta}\nu^{\beta}
\]
\\
\\
&=& \displaystyle 
- \[ 
\bar{\nabla}_{\alpha}(\nu^{\prime})^{\alpha} + {1\over 2}
\(
\bar{\nabla}_{\alpha} h^{\alpha}_{\beta} + 
\bar{\nabla}_{\beta} h^{\alpha}_{\alpha} -
\bar{\nabla}^{\alpha} h_{\alpha\beta}
\)\nu^{\beta}
\]
\\
\\
&=& \displaystyle 
- \bar{\nabla}_{\alpha}(\nu^{\prime})^{\alpha} - {1\over 2}
\bar{\nabla}_{\beta} (\Tr_{\bar{g}} h) \nu^{\beta}
\\
\\
&=& \displaystyle 
- \bar{\nabla}_{\alpha}(\nu^{\prime})^{\alpha} - {1\over 2}
\<\bar{\nabla}(\Tr_{\bar{g}}h),\nu\>_{\bar{g}}.
\end{array}
$$
Thus we obtain
\begin{equation}\label{eq4}
B_I=  \<\bar{\nabla}(\Tr_{\bar{g}}h),\nu\>_{\bar{g}} =-\( 2H^{\prime} + 2
 \bar{\nabla}_{\alpha}(\nu^{\prime})^{\alpha}\) .
\end{equation}
Now we compute the term $B_{II}$. We have
$$
\begin{array}{lcl}
\displaystyle 
\sum_{\alpha=0}^{n-1}
(\bar{\nabla}_{e_{\alpha}}h)(e_{\alpha}, \nu) &=&
\displaystyle 
\sum_{\alpha=0}^{n-1}\(
\bar{\nabla}_{e_{\alpha}} h(e_{\alpha},\nu) - 
h(\bar{\nabla}_{e_{\alpha}}e_{\alpha},\nu) -
h(e_{\alpha}, \bar{\nabla}_{e_{\alpha}}\nu)\)
\\
\\
&=&
\displaystyle 
\bar{\nabla}_{\nu} h(\nu,\nu) + \sum_{i=1}^{n-1}
\bar{\nabla}_{e_{i}} h(e_i,\nu) - h_{00}H + h^{ik} A_{ik}
\end{array}
$$
since 
$$
\begin{array}{l}
\displaystyle 
\sum_{\alpha=0}^{n-1}\bar{\nabla}_{e_{\alpha}}\nu= \bar{\nabla}_{\nu}\nu +
\sum_{i=1}^{n-1}\bar{\nabla}_{e_{i}}\nu, \ \ \ 
\bar{\nabla}_{\nu}\nu =0, \ \ \mbox{and}
\\
\\
\displaystyle 
\bar{\nabla}_{e_{i}}\nu= -\sum_{k=1}^{n-1}A_{ik}e_k, \ \ \ 
\sum_{\alpha=0}^{n-1}\bar{\nabla}_{e_{\alpha}}e_{\alpha} = 
\bar{\nabla}_{\nu}\nu + \sum_{i=1}^{n-1}\bar{\nabla}_{e_{i}}e_i=H\nu.
\end{array}
$$
Notice that $h^{ik}=fg^{ik}$ and $\bar{\nabla}_{e_{i}}  h(e_j,\nu) =
\nabla_{e_{i}}  h(e_j,\nu)$. Thus we have
\begin{equation}\label{eq5}
B_{II}= -\sum_{\alpha=0}^{n-1} (\bar{\nabla}_{e_{\alpha}}h)(e_{\alpha},
\nu) = - \bar{\nabla}_{\nu} h(\nu,\nu) - \sum_{i=1}^{n-1} \nabla_{e_{i}}
h(e_i,\nu) + h_{00}H - f H.
\end{equation}
To continue, we notice that $\bar{g}(\nu,\nu)=1$ implies
$$
0 = \bar{g}^{\prime}(\nu,\nu) + 2\bar{g}(\nu^{\prime},\nu) =
h(\nu,\nu) +  2\bar{g}(\nu^{\prime},\nu).
$$
Also we have
$$
0 = \bar{\nabla}_{\nu} h(\nu,\nu) + 
2 \bar{g}(\bar{\nabla}_{\nu} \nu^{\prime},\nu)  + 
2 \bar{g}(\nu^{\prime},\bar{\nabla}_{\nu} \nu)  = 
\bar{\nabla}_{\nu} h(\nu,\nu) + 
2 \bar{g}(\bar{\nabla}_{\nu} \nu^{\prime},\nu) 
$$
since $\bar{\nabla}_{\nu} \nu=0$. Thus we have that
\begin{equation}\label{eq6}
2\bar{g}(\bar{\nabla}_{\nu}\nu^{\prime},\nu) = 
-\bar{\nabla}_{\nu} h(\nu,\nu).
\end{equation}
Then the identity $\bar{g}(\nu,e_i)=0$ implies
$$
0 = \bar{g}^{\prime}(\nu,e_i) + \bar{g}(\nu^{\prime},e_i) + 
\bar{g}(\nu,e_i^{\prime}).
$$
Notice that $e_i^{\prime}\in T_pM$ since $g(t)\in C$. Thus
$$
0 = \bar{g}^{\prime}(\nu,e_i) + \bar{g}(\nu^{\prime},e_i).
$$
Now it follows that
$$
0 = \sum_{i=1}^{n-1}\(\bar{\nabla}_{e_i}h(\nu,e_i) +
\bar{g}(\bar{\nabla}_{e_i}\nu^{\prime},e_i) +
\bar{g}(\nu^{\prime},\bar{\nabla}_{e_i} e_i)\).
$$
Notice that $\sum_{i=1}^{n-1}\bar{\nabla}_{e_i} e_i=H\nu$, and
$H\bar{g}(\nu^{\prime},\nu)= -{1\over 2} h(\nu,\nu)$. Thus we obtain
\begin{equation}\label{eq7}
2\sum_{i=1}^{n-1}\bar{g}(\bar{\nabla}_{e_i}\nu^{\prime},e_i) = -2
\sum_{i=1}^{n-1}\nabla_{e_i}h(\nu,e_i) +  h_{00} H.
\end{equation}
We combine (\ref{eq6}) and (\ref{eq7}) to obtain
\begin{equation}\label{eq8}
\begin{array}{lcl}
2\bar{\nabla}_{\alpha}(\nu^{\prime})^{\alpha} &=&\displaystyle 2
\sum_{\alpha=0}^{n-1}\bar{g}(\bar{\nabla}_{e_{\alpha}}\nu^{\prime},e_{\alpha})
= 2 \[ \bar{g}(\bar{\nabla}_{\nu}\nu^{\prime}, \nu) +
\sum_{i=1}^{n-1}\bar{g}(\bar{\nabla}_{e_{i}}\nu^{\prime},e_{i}) \] 
\\ 
\\
&=&\displaystyle -\bar{\nabla}_{\nu} h(\nu,\nu) - 
2 \sum_{i=1}^{n-1}\nabla_{e_i}
h(\nu,e_i) + h_{00} H.
\end{array}
\end{equation}
Now it follows from (\ref{eq4}), (\ref{eq5}) and  (\ref{eq8}) that
$$
\begin{array}{lcl}
B_I + B_{II} &=& \displaystyle
-2H^{\prime} - 2 \bar{\nabla}_{\alpha}(\nu^{\prime})^{\alpha} -
\bar{\nabla}_{\nu}h(\nu,\nu) -
\sum_{i=1}^{n-1}\nabla_{e_i} h(e_i,\nu) + h_{00} H - f H
\\
\\
&=& \displaystyle 
-2H^{\prime} - \sum_{i=1}^{n-1}\nabla_{e_i} h(\nu,e_i)  - f H.
\end{array}
$$
Denote $\theta(v)= h(\nu,v) $ for $v\in T_xM$, so $\theta$ is a 1-form
on $M$. We notice that
$$
\nabla_{e_j}\theta(e_i) = (\nabla_{e_i}\theta)(e_j) +
\theta(\nabla_{e_i}e_j) = (\nabla_{e_i}\theta)(e_j) \ \ \ \mbox{since} \ \
\nabla_{e_i}e_j=0 \ \ \mbox{at $p$}.
$$
Thus we have that
$$
B_I + B_{II} = -2H^{\prime} - f H - \nabla_i\theta^i \ \ \ \mbox{on $M$}.
$$
This proves Claim \ref{Cl1} and concludes the proof of Theorem
\ref{Th1}. \hfill $\Box$
\section{Approximation Theorems}\label{AB-trick}
{\bf \ref{AB-trick}.1. Kobayashi approximation lemma.}  First we
reformulate several known facts in our terms. The following fact
follows from a modification of the continuity property of the Yamabe
constant due to B\'erard Bergery.  
\begin{Lemma}\label{app1} {\rm (cf. \cite[Proposition
4.31]{Besse})}
Let $\bar{g}_i$, $\bar{g}\in \Riem_C^0(W,M)$ be Riemannian metrics, and
$\bar{C}_i=[\bar{g}_i]$, $\bar{C}=[\bar{g}]$. Assume that
$$
\{
\begin{array}{ll}
\bar{g}_i \rightarrow \bar{g} & \mbox{in the $C^0$-topology on $W$, and}
\\
R_{\bar{g_i}} \rightarrow R_{\bar{g}} & \mbox{in the $C^0$-topology on $W$.}
\end{array}
\right.
$$
Then $Y_{\bar{C_i}}(W,M,C) \rightarrow Y_{\bar{C}}(W,M,C)$.
\end{Lemma}
Now we recall the results due to O. Kobayashi \cite{Kobayashi1}.
\begin{Lemma}\label{app2}{\rm (O. Kobayashi, \cite{Kobayashi1})}
For any $\delta>0$ there exists a smooth nonnegative function
$w_{\delta}$, and there exists $\epsilon(\delta)={1\over
4}e^{-{1\over\delta}}$ such that\\
\noindent
\parbox{4.0in}{
\begin{enumerate}
\item[{\bf (i)}] $
\{
\begin{array}{ll}
w_{\delta}(t)\equiv 1 & \mbox{on $[0,\epsilon(\delta)]$,}
\\
w_{\delta}(t)\equiv 0 & \mbox{on $[\delta,\infty)$},
\end{array}
\right.$
\item[{\bf (ii)}] \ \ $|t\dot{w}_{\delta}(t)| < \delta$ for $t\geq 0$, 
\item[{\bf (ii)}] \ \ $|t\ddot{w}_{\delta}(t)| < \delta$ for $t\geq 0$.
\end{enumerate}
{\rm (see Fig. \ref{AB-trick}.1.)}
}
\parbox{2.2in}{\PSbox{ryam-3.pstex}{10mm}{27mm}  
\begin{picture}(0,1)
\put(87,10){{\small $r$}}
\put(55,3){{\small $\delta$}}
\put(-10,3){{\small $\epsilon(\delta)$}}
\put(-22,45){{\small $1$}}
\put(-17,78){{\small $y$}}
\put(17,48){{\small $y=\omega_{\delta}(r)$}}
\end{picture}

\centerline{{\small {\bf Fig. \ref{AB-trick}.1}}}}
\end{Lemma}
\begin{Lemma}\label{app3}{\rm (O. Kobayashi, \cite{Kobayashi1})}
Let $\bar{g}, \tilde{g}\in \Riem(W)$, and $h= \tilde{g}-\bar{g}$. Then {\rm 
$$
\begin{array}{l}
R_{\tilde{g}} - R_{\bar{g}} = P_{\bar{g}}(h) + Q_{\bar{g}}(h), \ \
\mbox{where}
\\
\\
\{
\begin{array}{lcl}
P_{\bar{g}}(h) & = & -\Delta_{\bar{g}}(\Tr_{\bar{g}}h)+ \bar{\nabla}^i
 \bar{\nabla}^j h_{ij} - \<h,\Ric_{\bar{g}}\>_{\bar{g}}, 
\\
\\
 |Q_{\bar{g}}(h)| &\leq & C\( |\bar{\nabla}h|^2  q^3 + |h| \cdot
 |\bar{\nabla}^2 h| q^2 + (|h|\cdot |\bar{\nabla}^2 h| + |\Ric_{\bar{g}}|
\cdot |h|^2)q\),
\end{array}
\right.
\end{array}
$$
}where the constant $C>0$ depends only on $n=\dim W$, and $q\in
C^{\infty}_+(W)$ is a function satisfying $q \cdot \tilde{g}\geq
\bar{g}$.
\end{Lemma}
\begin{Proposition}\label{app4}
Let $W$ be a manifold with boundary $\p W= M$, and let metrics
$\bar{g}, \tilde{g}\in \Riem_C(W,M)$ such that $j^1_M\bar{g}=
j^1_M\tilde{g}$ (i.e. $\bar{g}$ coincides with $\tilde{g}$ up to
second derivatives on $M$), and $R_{\bar{g}}=R_{\tilde{g}}$ on
$M$. Then the family of metrics
$$
\tilde{g}_{\delta}= \bar{g} + w_{\delta}(r)(\tilde{g}-\bar{g})
\in \Riem_C(W,M) 
$$
satisfies the following properties:
\begin{enumerate}
\item[{\bf (i)}] $\tilde{g}_{\delta} \rightarrow \bar{g}$ in
the $C^1$-topology on $W$ {\rm (}as $\delta  \rightarrow 0${\rm )},
\item[{\bf (ii)}] $R_{\tilde{g}_{\delta}} \rightarrow R_{\bar{g}}$ in
the $C^0$-topology on $W$ {\rm (}as $\delta  \rightarrow 0${\rm )},
\item[{\bf (iii)}] $\tilde{g}_{\delta} \equiv \tilde{g}$ on the collar
$U_{\delta}(M,\bar{g})=\{ x\in W \ | \
\dist_{\bar{g}}(x,M)<\epsilon(\delta)\}$,
\item[{\bf (iv)}]  $\tilde{g}_{\delta} \equiv \bar{g}$ on $W\setminus 
U_{\epsilon(\delta)}(M,\bar{g})$.
\end{enumerate}
\end{Proposition}
\begin{Proof}
The statements (iii), (iv) are obvious. We prove (i) and (ii).

\noindent
{\bf (i)} The function $w_{\delta}$ is such that 
$\supp(w_{\delta})\subset [0,\delta]$. Then it follows
$$
\tilde{g}_{\delta} - \bar{g} = w_{\delta}(r)(\tilde{g} - \bar{g})
= O(r^2),
$$
thus $\tilde{g}_{\delta} \rightarrow \bar{g}$ in the $C^0$-topology on
$W$. Furthermore,
$$
\p(\tilde{g}_{\delta} - \bar{g}) = \dot{w}_{\delta}(r)
(\tilde{g}_{\delta} - \bar{g})+ w_{\delta}(r)\p (\tilde{g} -
\bar{g}).
$$
By the condition on the metrics $\tilde{g}_{\delta},\bar{g}$, 
$$
\tilde{g}_{\delta} - \bar{g} = O(r^2), \ \ \ \p(\tilde{g}_{\delta} -
\bar{g})= O(r).
$$
We use Lemma \ref{app2} to estimate
$$
|\p\tilde{g}_{\delta} - \p\bar{g}| \leq |\dot{w}_{\delta}(r)r|
\cdot { O(r^2)\over r} + w_{\delta}(r) \cdot O(r) \leq \delta
O(\delta) + O(\delta).
$$
Thus $\p\tilde{g}_{\delta} \rightarrow \p\bar{g}$ in the $C^0$-topology,
i.e. $\tilde{g}_{\delta} \rightarrow \bar{g}$ in the $C^1$-topology on $W$.

\noindent
{\bf (ii)} We use  Lemma \ref{app3} to write
$$
\{
\begin{array}{rcl}
R_{\tilde{g}_{\delta}}- R_{\bar{g}} &=& P_{\bar{g}}(w_{\delta}(r)
(\tilde{g}_{\delta} - \bar{g}))+ Q_{\bar{g}}(w_{\delta}(r)
(\tilde{g}_{\delta} - \bar{g})),
\\ 
\\
w_{\delta}(r)(R_{\tilde{g}_{\delta}}- R_{\bar{g}}) &=&
w_{\delta}(r)P_{\bar{g}} (\tilde{g}_{\delta} - \bar{g})+ 
w_{\delta}(r)Q_{\bar{g}}
(\tilde{g}_{\delta} - \bar{g}).
\end{array}
\right.
$$
We use again Lemma \ref{app3}:
$$
\begin{array}{l}
|P_{\bar{g}}(w_{\delta}(r) (\tilde{g}_{\delta} - \bar{g}))
- w_{\delta}(r)P_{\bar{g}} (\tilde{g}_{\delta} - \bar{g})|
\\
\\
\ \ \ \ \ \ \ \ 
\leq
C\(|\ddot w_{\delta}(r)|\cdot |w_{\delta}(r)| \cdot
|\tilde{g}_{\delta} - \bar{g}| + |\dot w_{\delta}(r)|^2\cdot
|\tilde{g}_{\delta} - \bar{g}|+
|\dot w_{\delta}(r)|\cdot |w_{\delta}(r)| \cdot
|\p(\tilde{g}_{\delta} - \bar{g})|\)
\\
\\
\ \ \ \ \ \ \ \ 
\leq
C\( \( |\ddot w_{\delta}(r)r^2| + 
|\dot w_{\delta}(r)r |^2 \) {O(r^2)\over r^2} + |\dot w_{\delta}(r)r |
\cdot {O(r)\over r}\) \leq C_1\delta.
\end{array}
$$
Similarly we obtain
$$
\{
\begin{array}{l}
|Q_{\bar{g}}(w_{\delta}(r)
(\tilde{g}_{\delta} - \bar{g}))| \leq C_2 \delta,
\\
\\
| w_{\delta}(r) Q_{\bar{g}}
(\tilde{g}_{\delta} - \bar{g}))| \leq C_3 \delta.
\end{array}
\right.
$$
Notice that 
$$
|w_{\delta}(r) (R_{\tilde{g}_{\delta}}- R_{\bar{g}})| \leq C_4\delta
$$
since $R_{\tilde{g}}\equiv R_{\bar{g}}$ on $M$. Thus we oblain:
$$
\begin{array}{lcl}
|R_{\tilde{g}_{\delta}}- R_{\bar{g}}| &\leq&
|P_{\bar{g}}(w_{\delta}(r) (\tilde{g}_{\delta} - \bar{g})) -
w_{\delta}(r)P_{\bar{g}} (\tilde{g}_{\delta} - \bar{g})| + \\ \\ 
& &
|Q_{\bar{g}}(w_{\delta}(r) (\tilde{g}_{\delta} - \bar{g}))| + 
|w_{\delta}(r) Q_{\bar{g}} (\tilde{g}_{\delta} - \bar{g}))| +
|w_{\delta}(r)(R_{\tilde{g}_{\delta}}- R_{\bar{g}})|\\ \\ 
&\leq&(C_1+C_2+C_3+C_4) \delta.
\end{array}
$$
Here $C_j$ ($j=1,\ldots,4$) are positive constants independent of
$\delta$. Thus $R_{\tilde{g}_{\delta}}\rightarrow R_{\bar{g}}$ in the
$C^0$-topology on $W$.
\end{Proof}
\begin{Proposition}\label{trick}
{\rm (Kobayashi Approximation Theorem, cf. \cite[Lemma 3.2]{Kobayashi1})} \\
Let $W$ be a manifold with boundary $\p W= M$, $C\in {\cal C}(M)$.
Let $\bar{g}\in \Riem_C(W,M)$ be a metric {\rm (}respectively $\bar{g}\in
\Riem_C^0(W,M)${\rm )}. Let $g= \bar{g}|_M$, and $A_{\bar{g}}$ be the second
fundamental form of $M=\p W$. There exists a family of metrics
$\tilde{g}_{\delta}\in \Riem_C(W,M)$ {\rm (}respectively 
$\tilde{g}_{\delta}\in
\Riem_C^0(W,M)${\rm )} such that
\begin{enumerate}
\item[{\bf (i)}] $\tilde{g}_{\delta} \rightarrow \bar{g}$ in
the $C^1$-topology on $W$ {\rm (}as $\delta  \rightarrow 0${\rm )},
\item[{\bf (ii)}] $R_{\tilde{g}_{\delta}} \rightarrow R_{\bar{g}}$ in
the $C^0$-topology on $W$ {\rm (}as $\delta  \rightarrow 0${\rm )},
\item[{\bf (iii)}] $\tilde{g}_{\delta}$ conformally equivalent to
the metric $(g-2rA_{\bar{g}}) + dr^2$ on 				
$U_{\epsilon(\delta)}(M,\bar{g})$,
\item[{\bf (iv)}]  $\tilde{g}_{\delta} \equiv \bar{g}$ on $W\setminus 
U_{\delta}(M,\bar{g})$.
\end{enumerate}
\end{Proposition}
\begin{Proof}
First, we note that the exponential map $\exp : T^{\perp} M
\longrightarrow W$ sends $(x,r\cdot\nu)\in T^{\perp} M $ to
$\exp_x(r\cdot\nu)=(x,r)\in W$. On $M$ we have
$$
\begin{array}{ll}
\bar{g}_{00}= \bar{g}(\p_r,\p_r) = 1, & \bar{g}_{00}^{\prime}= \p_r
\bar{g}(\p_r,\p_r) = 2 \bar{g}(\bar{\nabla}_{\p r}\p_r,\p_r)\equiv 0,
\\
\\
\bar{g}_{0i}= \bar{g}(\p_r,\p_i) = 0, & \bar{g}_{0i}^{\prime}= \p_r
\bar{g}(\p_r,\p_i) = \bar{g}(\bar{\nabla}_{\p r}\p_r,\p_i) + 
\bar{g}(\p_r,\bar{\nabla}_{\p r}\p_i) = 0 ,
\\
\\
\bar{g}_{ij}= \bar{g}(\p_i,\p_j) = g_{ij}, & \bar{g}_{ij}^{\prime}= \p_r
\bar{g}(\p_i,\p_j) = \bar{g}(\bar{\nabla}_{\p r}\p_i,\p_j) + 
\bar{g}(\p_i,\bar{\nabla}_{\p r}\p_j) = -2 A_{ij}.
\end{array}
$$

\noindent
\parbox{3.4in}{Here we used that
$$
\begin{array}{l}
\bar{\nabla}_{\p r}\p_r=0, \  \bar{\nabla}_{\p r}\p_i=-A_i^k\p_k,
\ \ \mbox{which implies}
\\
\\
\bar{g}(\p_r,\bar{\nabla}_{\p r}\p_i) = 
\bar{g}(\p_r,-A_i^k\p_k) =0, \ \ \mbox{and}
\\
\\
\bar{g}(\bar{\nabla}_{\p r}\p_i,\p_j) + 
\bar{g}(\p_i,\bar{\nabla}_{\p r}\p_j) 
\\
\\
= \bar{g}(\p_i,-A_j^k\p_k) +
\bar{g}(-A_i^k \p_k,\p_j) = -2 A_{ij}.
\end{array}
$$
}
\parbox{2.6in}{\hspace*{5mm}\PSbox{ryam-5.pstex}{10mm}{35mm}
\begin{picture}(0,1)
\put(-20,62){{\small $x$}}
\put(-5,67){{\small $\nu$}}
\end{picture}

\centerline{{\small {\bf Fig. \ref{AB-trick}.2.}}}}

\noindent
We define new metrics $\hat{g}$ and $G$ near $M$ and compare them with
the metric $\bar{g}$:
$$
\begin{array}{lcl}
\bar{g}(x,r) &=& 
(g_{ij}(x) - 2r A_{ij}(x) + O(r^2)) dx^idx^j + O(r^2) dr dx^i + dr^2,
\\
\\
\hat{g}(x,r) &:=& 
(g_{ij}(x) - 2r A_{ij}(x)) dx^idx^j +  dr^2,
\\
\\
G(x,r)  &:=& g_{ij}(x) dx^idx^j + dr^2.
\end{array}
$$
Clearly $j^1_M \bar{g} = j^1_M \hat{g}$, and, in general, 
$j^1_M \bar{g} \neq j^1_M G$. We notice 
$$
\begin{array}{ccl}
R_{\bar{g}}|_M & = & 
R_G + 2\Ric_{\bar{g}}(\nu,\nu) + |A_{\bar{g}}|^2_g - H^2_{\bar{g}}, \ \ 
\\
\\
& = & 
R_g + 2\Ric_{\bar{g}}(\nu,\nu) + |A_{\bar{g}}|^2_g - H^2_{\bar{g}}.
\end{array}
$$
We define a metric 
$$
\check{g}= \(\check{g}_{ij}\):= \(g_{ij}(x) -2rA_{ij}(x)\)
$$
on each hyperseface $M\times \{r\}\subset W$ (for small $r$). 
Then we have
$$
\begin{array}{lcl}
R_{\hat{g}} &=& \displaystyle
R_{\check{g}} + {3\over 4}|\p_r
\check{g}_{ij}|^2_{\check{g}} - \check{g}^{ij}\cdot \p_r^2 \check{g}_{ij}
- {1\over 4} | \check{g}_{ij}\cdot \p_r \check{g}_{ij}|^2
\\
\\
&=& \displaystyle
R_g+ 3|A_{\bar{g}}|_g^2 - H_{\bar{g}}^2 + O(r) \ \ \ \mbox{near $M$, and}
\\
\\
R_{\hat{g}} &=& \displaystyle
R_g+ 3|A_{\bar{g}}|_g^2 - H_{\bar{g}}^2 \ \ \ \mbox{on $M$.}
\end{array}
$$
We choose the conformal metric $\tilde{g}(x,r)= u(x,r)^{{4\over n-2}}
\cdot \hat{g}$ such that $j^1_M\bar{g}=j^1_M\tilde{g}$ by giving $u$ the
boundary conditions:
$$
u(x,0)\equiv 1, \ \  \ \p_ru(x,0)\equiv 0 \ \ \mbox{on $M$}.
$$
We have 
$$
\begin{array}{l}
\displaystyle
-{4(n-1)\over n-2} \Delta_{\hat{g}} u + R_{\hat{g}}u =
 R_{\tilde{g}}u^{{n+2\over n-2}}, \ \ \mbox{or}
\\
\\
\displaystyle 
\Delta_{\hat{g}} u = -{n-2\over 4(n-1)}\( 
R_{\tilde{g}}u^{{n+2\over n-2}} -  R_{\hat{g}}u\).
\end{array}
$$
We specify $\Delta_{\hat{g}} u$ on $M$:
$$
\begin{array}{lcl}
\displaystyle
\Delta_{\hat{g}} u &=& \displaystyle
\nabla^{\alpha}\p_{\alpha} u = \hat{g}^{\alpha\beta}\(\p_{\alpha}\p_{\beta}u
-\hat{\Gamma}^{\gamma}_{\alpha\beta}\p_{\gamma}u\)
\\
\\
&=& \displaystyle
\p_r^2 u + g^{ij}\p_i\p_j u - \hat{\Gamma}_{00}^{0}\p_r u - 
\hat{\Gamma}_{00}^{i}\p_i u - g^{ij}\(\hat{\Gamma}_{ij}^{0}\p_r u - 
\hat{\Gamma}_{ij}^{k}\p_k u\) = \p_r^2 u
\end{array}
$$
since $\p_i\p_j u =0$, $\p_r u=0$, and $\p_iu=0$ on $M$. Here we use the
$u(x,0)\equiv 1$, $\p_ru(x,0)\equiv 0$ on $M$. Thus we obtain that on
$M$
$$
\begin{array}{lcl}
\p^2_r u &=& \displaystyle \Delta_{\hat{g}} u = -{n-2\over 4(n-1)}\( 
R_{\tilde{g}} -  R_{\hat{g}} \) 
\\
\\
&=& \displaystyle 
-{n-2\over 4(n-1)}\( 
R_{\tilde{g}} -  (R_g + 3|A_{\bar{g}}|_g^2 -H_g^2)\).
\end{array}
$$
We let $u(x,r):= 1 + {1\over 2} r^2 \phi(x)$ near $M$, where
\begin{equation}\label{a1}
\begin{array}{rcl}
\phi(x) &= & \displaystyle
-{n-2\over 4(n-1)}\( 
R_{\bar{g}}|_M -  (R_g + 3|A_{\bar{g}}|_g^2 -H_g^2)\)
\\
\\
&= & \displaystyle
-{n-2\over 4(n-1)}\( 
R_{\bar{g}}|_M -  R_{\hat{g}} \) .
\end{array}
\end{equation}
Then the metric 
$$
\tilde{g}= u^{{4\over n-2}} \cdot \hat{g} = \(1 + {1\over 2} r^2
\phi (x)\)^{{4\over n-2}}\[ (g -2rA) + dr^2\]
$$
is such that $j^1_M\tilde{g}=j^1_M\bar{g}$, and
$R_{\tilde{g}}=R_{\bar{g}}$ on $M$. We use Proposition \ref{app4} to define
a family of metrics $\tilde{g}_{\delta}$:
$$
\tilde{g}_{\delta} = \bar{g} + w_{\delta}(r) \cdot (\tilde{g}-\bar{g}) \in
\Riem_C(W,M).
$$
We also notice that 
$$
\begin{array}{lcl}
C= [\bar{g}|_M]=
[\hat{g}|_M]=[\tilde{g}|_M]=[\tilde{g}_{\delta}|_M], \ \ \mbox{and}
\\
\\
H_{\bar{g}} = 0 \Longrightarrow H_{\hat{g}} = 0 
\Longrightarrow H_{\tilde{g}} = 0 \Longrightarrow H_{\tilde{g}_{\delta}} = 0 
\end{array}
$$
since $A_{\bar{g}}= A_{\hat{g}}$, $\p_r u= 0$ on $M$, and
$\tilde{g}=\tilde{g}_{\delta}$ near $M$. Then $\bar{g}\in
\Riem_C^0(W,M)$ implies that $\tilde{g}_{\delta}\in \Riem_C^0(W,M)$.
\end{Proof} 
{\bf \ref{AB-trick}.2. The approximation trick under minimal boundary
condition.} One notices that the above results do not allow to use a
metric which is conformally equivalent to a product metric near the
boundary to approximate the relative Yamabe constant
$Y_{\bar{C}}(W,M;C)$. This is the problem which we address and solve
here.
\begin{Theorem}\label{main-trick}
{\rm (Approximation Trick)} \\ Let $W$ be a manifold with boundary $\p
W= M$, $C\in {\cal C}(M)$.  Let $\bar{g}\in \Riem_C^0(W,M)$ be a
metric. Let $g= \bar{g}|_M$, and $A_{\bar{g}}$ be the second
fundamental form of $M=\p W$. There exists a family of metrics
$\tilde{g}_{\delta}\in \Riem_C^0(W,M)$ such that
\begin{enumerate}
\item[{\bf (i)}] $\tilde{g}_{\delta} \rightarrow \bar{g}$ in
the $C^0$-topology on $W$ {\rm (}as $\delta  \rightarrow 0${\rm )},
\item[{\bf (ii)}] $R_{\tilde{g}_{\delta}} \rightarrow R_{\bar{g}}$ in
the $C^0$-topology on $W$ {\rm (}as $\delta  \rightarrow 0${\rm )},
\item[{\bf (iii)}] $\tilde{g}_{\delta}$ conformally equivalent to
the metric $g + dr^2$ on 				
$U_{\epsilon(\delta)}(M,\bar{g})$,
\item[{\bf (iv)}]  $\tilde{g}_{\delta} \equiv \bar{g}$ on $W\setminus 
U_{\delta}(M,\bar{g})$.
\end{enumerate}
\end{Theorem}
{\bf Remark.} In order to control the scalar curvature without the
minimal boundary condition, one needs the $C^1$-convergence of metrics
as in Proposition \ref{trick}. Furthemore, when $\bar{g}$ is not
totally umbilic on $M$, the metric $\bar{g}$ can never be approximated
in the $C^1$-topology to a metric which conformally is a product
metric near the boundary. However, we emphasize that the convergence
in (i) of Theorem \ref{main-trick} is the $C^0$-convergence only. The
minimal boundary condition plays a crucial role to achieve the
$C^0$-convergence for scalar curvatures in (ii).
\begin{Proof} There are two steps in the proof.
\vspace{2mm}

\noindent
{\bf Step 1.} First, Proposition \ref{trick} allows us to assume that
the metric $\bar{g}$ is such that
$$
\bar{g} = \(1+ {r^2\over 2} \phi (x)\)^{{n-2\over 4}} \[ (g(x) -2 r
A_{\bar{g}}(x)) +dr^2\] \ \ \ \mbox{on a collar
$U_{\delta_0}(M,\bar{g})$},
$$
where $\{x,r\}=\{x^1,\ldots,x^{n-1},r\}$ denotes a Fermi coordinate
system near each point of $M$, and $\phi(x)$ is the $C^{\infty}$-function
on $M$ defined by (\ref{a1}).
\vspace{2mm}

\noindent
For each positive $\delta <\delta_0$, let $G_{\delta}\in
\Riem^0_C(W,M)$ be a metric defined by
$$
\begin{array}{rcl}
G_{\delta}(x,r) &=& \hat{g}(x,r) + w_{\delta}(r) \cdot (G(x,r) -
\hat{g}(x,r)) \\ \\ &=& g(x) - 2r(1-w_{\delta}(r) )\cdot A_{\bar{g}} +
dr^2.
\end{array}
$$
Here $\hat{g}(x,r)$ and $G(x,r)$ are given by
$$
\left.
\begin{array}{rcl}
\hat{g}(x,r)&=& 
(g(x) - 2rA_{\bar{g}}(x)) + dr^2,
\\
\\
G(x,r)&=& 
g(x) + dr^2 \ 
\end{array}
\} \ \ \mbox{on $U_{\delta}(M,\bar{g})$}.
$$
We also let $\check{g}_{\delta}(x,r) = g(x) -2r(1-w_{\delta}(r))\cdot
A_{\bar{g}}(x)$ on $U_{\delta}(M,\bar{g})$. It follows then from Lemma
\ref{app2} that near $M$ the scalar curvature of the metric
$G_{\delta}$ satisfies
$$
\begin{array}{rcl}
R_{G_{\delta}} & = & \displaystyle
R_{\check{g}} + {3\over 4}|\p_r
(\check{g}_{\delta})_{ij}|^2_{\check{g}_{\delta}} -
\check{g}_{\delta}^{ij}\cdot \p_r^2(\check{g}_{\delta})_{ij}
-{1\over 4}|\check{g}_{\delta}^{ij}\cdot \p_r (\check{g}_{\delta})_{ij}|^2
\\
\\
& = & \displaystyle
R_g + 3(1 -w_{\delta}(t))^2|A_{\bar{g}}|^2_g
-(1 -w_{\delta}(t))^2H_{\bar{g}}^2 - (4\dot{w}_{\delta}(r) + 
2r \cdot \ddot{w}_{\delta}(r))H_{\bar{g}} + O(\delta).
\end{array}
$$
We use the minimal boundary condition $H_{\bar{g}}=0$ to obtain
$$
R_{G_{\delta}}= R_g + 3(1-w_{\delta}(r))^2 |A_{\bar{g}}|^2_g + O(\delta)
\ \ \ \ \mbox{near $M$.}
$$
{\bf Step 2.} We define now the metric $\tilde{g}_{\delta}\in
\Riem_C^0(W,M)$ as follows:
\begin{equation}\label{a2}
\tilde{g}_{\delta}(x,r)= \(1+{r^2\over 2} \phi_{\delta}(x,r) \)^{{4\over
n-2}} \cdot G_{\delta}(x,r)
\end{equation}
on $U_{\delta}(M,\bar{g})$, with
\begin{equation}\label{a3}
\phi_{\delta}(x,r) = \phi(x) - {3(n-2)\over 4(n-1)}
(2-w_{\delta}(r)) w_{\delta}(r) |A_{\bar{g}}|^2_g .
\end{equation}
We obtain that the assertions (iii) and (iv) hold since $G_{\delta}=
g+dr^2$ on the collar $U_{\epsilon(\delta)}(M,\bar{g})$, and
$G_{\delta}= \hat{g}$, and $\phi_{\delta}=\phi$ outside of the collar
$U_{\delta}(M,\bar{g})$. By construction
$$
G_{\delta}\lra \hat{g}, \ \ \ \mbox{and} \ \ \ {r^2\over
2}\phi_{\delta}(x,r) \lra {r^2\over
2}\phi(x)
$$
in the $C^0$-topology on $W$ as $\delta\lra 0$. Thus the assertion (i) holds.
Finally, the scalar curvature $R_{\tilde{g}_{\delta}}$ is given by
$$
\begin{array}{rcl}
R_{\tilde{g}_{\delta}}
& = & 
(1+{r^2\over 2} \phi_{\delta}(x,r))^{-{n+2\over n-2}}
\[
-{4(n-2)\over n-1}\Delta_{G_{\delta}} (1+{r^2\over 2} \phi_{\delta}(x,r))
+ R_{G_{\delta}} (1+{r^2\over 2} \phi_{\delta}(x,r))
\]
\\
\\
& = & 
(1+ O(\delta^2))
\[
-{4(n-2)\over n-1} \phi + 3(2-w_{\delta}(r))\cdot w_{\delta}(r)\cdot
|A_{\bar{g}}|^2_g +
R_{G_{\delta}} + O(\delta)
\]
\\
\\
& = & 
R_{G_{\delta}} + (R_{\bar{g}} -R_g - 3 |A_{\bar{g}}|^2_g) + 3(2-w_{\delta}(r))\cdot w_{\delta}(r)\cdot
|A_{\bar{g}}|^2_g+ O(\delta)
\\
\\
& = & 
R_g + 3(1-w_{\delta}(r))^2 |A_{\bar{g}}|_g^2 + 
(R_{\bar{g}} -R_g - 3 |A_{\bar{g}}|^2_g) + 
3(2-w_{\delta}(r)) w_{\delta}(r)
|A_{\bar{g}}|^2_g+ O(\delta)
\\
\\
& = & 
R_{\bar{g}} + O(\delta) \ \ \ \mbox{on $W$ as $\delta\lra 0$}.
\end{array}
$$
This implies the assertion (ii) and completes the proof of Theorem
\ref{main-trick}.
\end{Proof}
\section{Gluing Theorems}\label{diss}
{\bf \ref{diss}.1. Setting.} Here we would like to analyze the gluing
procedure for manifolds equipped with conformal structures.  Let
$W_1$, $W_2$ be two compact manifolds, $\dim W_1 = \dim W_2 \geq 3$,
with boundaries
$$
\p W_1 = M_1= M_0 \sqcup M,  \ \ \ \mbox{and}\ \ \ 
\p W_2 = M_2= M_0 \sqcup M^{\prime}
$$
endowed with conformal classes $C_1 = C_0\sqcup C\in {\cal C}(M_1)$,
$C_2 = C_0\sqcup C^{\prime}\in {\cal C}(M_2)$, where $C_0\in {\cal
C}(M_0)$, $C\in {\cal C}(M)$, $C^{\prime}\in {\cal C}(M^{\prime})$.
Let $W= W_1\cup_{M_0}(-W_2)$ be the boundary connected sum 
of $W_1$ and $W_2$ along $M_0$ (see Fig. \ref{diss}.1).
\vspace{2mm}

\noindent
{\bf Remark.} The boundary of the manifold $W$ is $\p W = M \sqcup
M^{\prime}$ with appropriate orientation. We consider both cases when
$\p W =\emptyset$ and $\p W \neq\emptyset$. 
\vspace{2mm}

\noindent
Recall that a conformal class $C\in {\cal C}(M)$ is {\it positive} if
$Y_C(M)>0$. 

\hspace*{20mm}\PSbox{ryam-1.pstex}{10mm}{64mm}
\begin{picture}(0,0)
\put(120,160){$M_0$}
\put(85,160){$M_0$}
\put(30,155){$W_1$}
\put(170,155){$W_2$}
\put(-45,130){$M$}
\put(285,125){$M^{\prime}$}
\put(-25,30){$M$}
\put(270,25){$M^{\prime}$}
\put(45,55){$W_1$}
\put(155,55){$W_2$}
\put(105,63){$M_0$}
\end{picture}

\centerline{{\small {\bf Fig. \ref{diss}.1.} Manifold $W=
W_1\cup_{M_0}(-W_2)$.}}
\begin{Theorem}\label{ThI}
Let $C_0\in {\cal C}(M_0)$ be a positive conformal class, and
$Y(W_j,M_j;C_j)> 0 $ for $j=1,2$. Then $Y(W,\p W; C\sqcup
C^{\prime})>0$.
\end{Theorem}
{\bf  Remark.} We do not assume that the conformal classes $C\in {\cal
C}(M)$, $C^{\prime}\in {\cal C}(M^{\prime})$ are positive.
\begin{pf-of}{\ref{diss}.2. Proof of Theorem \ref{ThI}.} 
There are four steps in the proof.
\vspace{2mm}

\noindent
{\bf Step 1.} First we notice that since $C_0\in {\cal C}(M_0)$ is a
positive conformal class, there exists a metric $h\in C_0$ on $M_0$
with $R_h>0$. The metric $h$ do not have to be a Yamabe metric. We fix
the metric $h$. The condition $Y(W_j,M_j;C_j)> 0 $ (for $j=1,2$)
implies that there exist conformal classes $\bar{C}_j$ on $W_j$ so
that $\p \bar{C}_j = C_j$, i.e. $(\bar{C}_j,C_j)\in {\cal C}(W_j,M_j)$.
We denote
$$
Y_{\bar{C}_j}= Y_{\bar{C}_j}(W_j,M_j;C_j)> 0 ,  \ \ \ j=1,2.
$$
Let $\bar{g}_j\in \bar{C}_j$ be such that
$\bar{g}_j|_{M_0}=h$. Moreover, we may assume that $\bar{g}_j\in
\bar{C}_j^{0}$ (i.e. that $H_{\bar{g}_j}\equiv 0$ on $M_j$). 
\vspace{2mm}

\noindent
{\bf Remarks. (1)} The metrics $\bar{g}_j\in \bar{C}_j^{0}$ do not have to
be relative Yamabe metrics, moreover, their scalar curvature $R_{\bar{g}_j}$ is
not positvive, in general.
 
\noindent
{\bf (2)} The union $\bar{C}_1^0\cup_C\bar{C}_2^0$ does not make sense
as a conformal class on $W$ since this union, in general, fails to be
smooth along $M_0$.
\vspace{2mm}

\noindent
{\bf Step 2.}  Theorem \ref{main-trick} and (\ref{a2}), (\ref{a3})
imply that for any $\epsilon>0$ there exist conformal classes
$\hat{C}_j$ on $W_j$ and metrics $\hat{g}_j\in \hat{C}_j$ ($j=1,2$)
such that
$$
\{
\begin{array}{l}
\ \p \hat{C}_j = C_j,
\\
\left.
\begin{array}{l}
\hat{g}_j \sim \bar{g}_j,
\\
R_{\hat{g}_j} \sim R_{\bar{g}_j}, 
\end{array}\} \ \ C^0\mbox{-close on} \ W_j, \  \
\mbox{which implies $|Y_{\hat{C}_j}-
Y_{\bar{C}_j}| < \epsilon$.}
\\
\ \hat{g}_j =\(1+ {r^2\over 2} f_j\)^{4\over n-2} \cdot (h+dr^2) \ \ 
\mbox{near $M_j$ in $W_j$.}
\end{array}
\right.
$$
Here the function $f_j$ is defined by
$$
f_j = -\frac{n-2}{4(n-1)}\(R_{\bar{g}_j}|_{M_0} - R_{h} + 
3|A_{\bar{g}_j}|^2_{h}\) \ \ 
\mbox{on $M_0$ in each $W_j$}.
$$
From now on we only need the conditions $Y_{\hat{C}_j}>0$ ($j=1,2$).
Therefore we may assume that $f_j<0$ on $M_0$ since the relative
Yamabe constant $Y_{\hat{C}_j}$ is invariant under pointwise
conformal change. 
\vspace{2mm}

\noindent
\hspace*{10mm}
\PSbox{ryam-2.pstex}{10mm}{30mm}
\begin{picture}(0,0)
\put(-40,30){$M$}
\put(360,25){$M^{\prime}$}
\put(45,55){$W_1$}
\put(230,55){$W_2$}
\put(125,47){$M_0\times [0,\ell]$}
\end{picture}
\vspace{2mm}

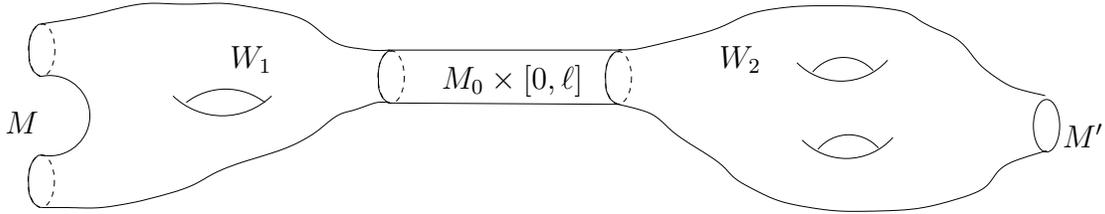
\centerline{{\small {\bf Fig. \ref{diss}.2.}  Manifold $X=
W_1\cup_{M_0}\( M_0\times [0,\ell]\) \cup_{M_0} (-W_2)$.}}
\vspace{2mm}

\noindent
Let $\ell$ be a positive constant. We define the
manifold $X$ which is diffeomorphic to $W$ as follows (see
Fig. \ref{diss}.2):
$$
X= W_1\cup_{M_0}\( M_0\times [0,\ell]\) \cup_{M_0} (-W_2)
$$
Now we need the cut-off function $w _{\delta}$ defined in Lemma
\ref{app2}. Then for each $\delta$, $0< \delta< \ell$, we define a
metric $\tilde{g}$ on $X$ as follows.
$$
\tilde{g} = 
\{
\begin{array}{ll}
\hat{g}_1 \ \ &\mbox{on $W_1$,}
\\
\hat{g}_2 \ \ &\mbox{on $W_2$,}
\\
(1+ {r^2\over 2} w _{\delta}(r)f_1)^{4\over n-2}(h+dr^2) 
\ \ &\mbox{on $M_0\times [0,\delta]$,}
\\
h+ dr^2 \ \ &\mbox{on $M_0\times [\delta,\ell-\delta]$,}
\\
(1+ {(\ell-r)^2\over 2} w _{\delta}(\ell-r)f_2)^{4\over n-2}(h+dr^2) 
\ \ &\mbox{on $M_0\times [\ell-\delta,\ell]$.}
\end{array}
\right.
$$
Clearly $\tilde{g}$ is a $C^{\infty}$-metric on $X\cong W$. Let
$\tilde{C}=[\tilde{g}]\in {\cal C}(W)$.
\vspace{2mm}

\noindent
{\bf Remark.} The metric $\tilde{g}$ does not have, in general,
positive scalar curvature.
\vspace{2mm}

\noindent
{\bf Step 3.} Let $j=1,2$. Denote $\nu_j$ the first eigenvalue of the
Yamabe operator on $W_j$ for the Neumann boundary condition. Then
\vspace{-.5em}
$$
\begin{array}{c}\displaystyle
\nu_j = \inf_{\begin{array}{c} _{u\in C^{\infty}(W_j)} \\
^{u\not\equiv 0}\end{array}} {\int_{W_j} \( {4(n-1)\over n-2} |d
u|^2_{\hat{g}_j} + R_{\hat{g}_j} u^2 \) d\sigma_{\hat{g}_j}\over
\int_{W_j} u^2 d\sigma_{\hat{g}_j} } \ .
\end{array}
$$
The relative Yamabe constants $Y_{\hat{C}_j}>0$ since $Y_{\hat{C}_j} >
Y_{\bar{C}_j} -\epsilon$, $j=1,2$. Thus it follows that $\nu_j >0$.
Notice that the conditions $R_h >0$ on $M$ and $f_j<0$ ($j=1,2$) imply
that $R_{\tilde{g}}>0$ on the cylinder $M_0\times [0,\ell]$ for small
$\delta>0$.
\vspace{2mm}

\noindent
Let $\nu_{cyl}$ be the first eigenvalue of the Yamabe operator on
$M_0\times [0,\ell]$ for the Neumann boundary condition. We have
\vspace{-.5em}
$$
\nu_{cyl} = \inf_{\begin{array}{c} _{u\in C^{\infty}(M_0\times [0,\ell])}
\\ ^{u\not\equiv 0}\end{array}} {\int_{M_0\times [0,\ell]} \(
{4(n-1)\over n-2} |d u|^2_{\tilde{g}} + R_{\tilde{g}} u^2 \)
d\sigma_{\tilde{g}}\over \int_{M_0\times [0,\ell]} u^2
d\sigma_{\tilde{g}} } \ .
$$
It follows that $\nu_{cyl}>0$ since $R_{\tilde{g}}>0$.
\vspace{2mm}

\noindent
{\bf Step 4.} Let $\nu$ be the first eigenvalue of the Yamabe operator on
$X\cong W$ for the Neumann boundary condition, which is equal to
$$
\nu = \inf_{\begin{array}{c} _{u\in C^{\infty}(X)}
\\ ^{u\not\equiv 0}\end{array}} {\int_{X} \(
{4(n-1)\over n-2} |d u|^2_{\tilde{g}} + R_{\tilde{g}} u^2 \)
d\sigma_{\tilde{g}}\over \int_{X} u^2
d\sigma_{\tilde{g}} }\ .
$$
We conclude that $\nu\geq
\min\{\nu_1,\nu_2,\nu_{cyl}\}>0$ by \cite[pp. 18-19]{Chavel}. The
condition $\nu>0$ is equivalent that there exists a metric
$\tilde{\tilde{g}}\in [\tilde{g}]$ such that $R_{\tilde{\tilde{g}}}>0$
on $X\cong W$ and $H_{\tilde{\tilde{g}}}=0$ on $\p X = \p W$. Thus $
Y_{[\tilde{\tilde{g}}]}(X,\p X; C\sqcup C^{\prime}) > 0$, and this
implies $Y (W,M\sqcup M^{\prime}; C\sqcup C^{\prime}) > 0 $.
\end{pf-of}

\noindent 
{\bf Remark.} We emphasize that we can choose $\ell > 0$ to be
small by choosing small $\delta > 0$. 
\vspace{2mm}

\noindent
{\bf \ref{diss}.3. Manifolds with positive Yamabe invariant.}
Here we would like to show that there are many examples of manifolds
with positive relative Yamabe invariant. 
\vspace{2mm}

\noindent
We start with a closed compact manifold $N$, $\dim N\geq 3$ with
$Y(N)>0$. We choose a small disk $D^n\subset N$ centered at $x_0\in
N$, then $\p(N\setminus \iint (D^n))= S^{n-1}$. Let $C_0\in {\cal
C}(S^{n-1})$ be the standard conformal class. (In particular, $C_0$ is
a positive class.) 
\begin{Theorem}\label{ThII}
Let $N$ be a closed compact manifold, $\dim N\geq 3$ with $Y(N)>0$.
Then $Y(N\setminus \iint (D^n), S^{n-1};C_0)>0$.
\end{Theorem}
\begin{Proof}  We use
\cite[Corollary 3.5.]{Kobayashi1} to choose a conformal class
$\bar{C}\in {\cal C}(N)$ with the Yamabe constant $Y_{\bar{C}}(N)>0$,
and a metric $\bar{g}\in \bar{C}$, such that

$\bullet$ $\displaystyle \bar{g}$ is conformally flat near $x_0\in N$,

$\bullet$ $R_{\bar{g}}>0$ on $N$.

\noindent
Thus (as it was observed, say, by Gromov-Lawson \cite{GL1}), there
exists a metric $\hat{g}$ on the manifold $N\setminus \iint (D^n)$ such
that
\begin{enumerate}
\item[$\bullet$] $\p[\hat{g}]= C_0\in {\cal C}(S^{n-1})$,
\item[$\bullet$] $R_{\hat{g}}>0$ on $N\setminus \iint (D^n)$,
\item[$\bullet$] $\displaystyle \hat{g} = g_{S^{n-1}}+dr^2$ near
$S^{n-1}=\p (N\setminus \iint (D^n))$, $[g_{S^{n-1}}]=C_0$.
\end{enumerate}
Thus $Y_{[\hat{g}]} (N\setminus \iint (D^n), S^{n-1};C_0)>0$ and
$Y(N\setminus \iint (D^n), S^{n-1};C_0)>0$.
\end{Proof}
{\bf \ref{diss}.4. The double.}  Let $W$ be a compact manifold with
$\p W = M= M_0 \sqcup M_1$. We define the manifold $X=
W\cup_{M_0}(-W)$, the double of $W$ along $M_0$ (see
Fig. \ref{diss}.3).
\vspace{2mm}

\noindent
{\bf Remark.} Here the other boundary component $M_1$ of $\p W$ may be
empty or not. 
\vspace{2mm}

\PSbox{ryam-4.pstex}{10mm}{25mm}  
\begin{picture}(0,1)
\put(-43,40){$M_1$}
\put(87,53){$M_0$}
\put(45,38){$W$}
\put(205,38){$W$}
\put(265,38){$-W$}
\put(235,55){$M_0$}
\end{picture}
\vspace{2mm}

\centerline{{\small {\bf Fig. \ref{diss}.3.}}}
\vspace{2mm}

\noindent
Theorem \ref{ThI} immediately implies the following.
\begin{Corollary}\label{ThIII}
Let $C=C_0\sqcup C_1\in {\cal C}(\p W)$ $( C_0\in {\cal C}(M_0) )$ and
$C_1\in {\cal C}(M_1)$. Let $Y(W,\p W; C) >0$. Then $Y(X,M_1 \sqcup
(-M_1); C_1\sqcup C_1) >0$.
\end{Corollary}
\section{Non-positive Yamabe invariant}\label{non-pos}
{\bf \ref{non-pos}.1. Setting.}  Let $W$ be a compact smooth
$n$-manifold with boundary $\p W= M\neq \emptyset$, and $M= M_{0}
\sqcup M_{1}$.  (Here $M_{1}$ may be empty.) Let $(\bar{C}, C = C_{0}
\sqcup C_{1})$ be a pair of conformal classes on $(W, M)$, where
$C_0\in {\cal C}(M_0)$, $C_1\in {\cal C}(M_1)$. Similar to the case of
closed manifolds, we first notice the following.
\begin{Proposition}\label{[Proposition 4.12]}{\rm (cf. 
\cite[Lemma 1.6]{Kobayashi1})}
Suppose $Y_{\bar{C}}(W, M; C)\leq 0$.  Then, for any $\bar{g}\in
\bar{C}^{0}$,
$$ 
(\min R_{\bar{g}}) \Vol_{\bar{g}}(W)^{2/n}  \leq   
Y_{\bar{C}}(W, M; C)  \leq     
(\max R_{\bar{g}})\Vol_{\bar{g}}(W)^{2/n}.  
$$
\end{Proposition}
Let $X = W \cup_{M_{0}}(-W)$ be the double of $W$ along $M_{0}$.
Then, the boundary of $X$ is $\p X= M_{1} \sqcup (-M_{1})$.  We use
Proposition \ref{[Proposition 4.12]} and Theorem \ref{main-trick} to
prove the following result.
\begin{Theorem}\label{[Theorem 4.12]}$\mbox{\ }$
\begin{enumerate}
\item[{\bf (1)}] If  $Y_{\bar{C}}(W, M; C)\leq 0$, then

\hspace*{10mm}$\displaystyle
2^{{2\over n}} Y_{\bar{C}}(W, M; C)  \leq  
    Y(X, M_{1} \sqcup (-M_{1}); C_{1}  \sqcup  C_{1}).
$
\item[{\bf (2)}] If  $Y_{\bar{C}}(W, M; C)\leq 0$, then

\hspace*{10mm}$\displaystyle
    2^{{2\over n}} Y(W, M; C)  \leq 
    Y(X, M_{1} \sqcup (-M_{1}); C_{1}  \sqcup C_{1}).  
$
\item[{\bf (3)}]   If  $Y(W, M)\leq 0$, then

\hspace*{10mm}$  \displaystyle 
2^{{2\over n}} Y(W, M)   \leq   Y(X, M_{1} \sqcup (-M_{1})). 
$
\end{enumerate}
\end{Theorem}
\begin{Proof}
Clearly $(1) \Longrightarrow (2) \Longrightarrow (3)$. Thus it is
enough to prove (1). Notice that if $Y(X, M_{1} \sqcup (-M_{1}); C_{1}
\sqcup C_{1})>0$ then (1) holds trivially.
\vspace{2mm}

\noindent
Let $Y(X, M_{1} \sqcup (-M_{1}); C_{1} \sqcup C_{1})\leq 0$. We choose
a conformal class $\bar{C}\in {\cal C}(W)$.  Then for a generic class
$\bar{C}$ there is a relative Yamabe metric $\bar{g}\in
\bar{C}^0$. Theorem \ref{main-trick} implies that, for any
$\epsilon>0$, there exists a metric $\hat{g}\in \hat{C}$ such that $g
:=\bar{g}|_M= \hat{g}|_M$, and
$$
\{
\begin{array}{l}
\bar{g} \sim \hat{g} \ \ C^0\mbox{-close on $W$, which implies 
$|\Vol_{\hat{g}}(W)^{{2\over n}} -\Vol_{\bar{g}}(W)^{{2\over n}}|<\epsilon$,} 
\\ 
|R_{\bar{g}} -
R_{\hat{g}}| <\epsilon  \ \ \mbox{on $W$,}
\\
\hat{g} = (1+{r^2\over 2}f(x))^{{4\over n-2}}(g+ dr^2) \ \ 
\mbox{near $M\subset W$}.
\end{array}
\right.
$$
We define $\tilde{g}:= \hat{g}\cup \hat{g}$ on $X= W\cup_{M_0}(-W)$.
The metric $\tilde{g}$ is smooth by construction. Let $\tilde{C}:=
\hat{C}\cup \hat{C}\in {\cal C}(X)$. By Proposition \ref{[Proposition
4.12]} we have
$$
\begin{array}{lcl}
Y_{\tilde{C}}(X,M_1\sqcup
(-M_1); C_1\sqcup C_1) &\geq & \displaystyle
\(\min_{X} R_{\tilde{g}}\)
\Vol_{\tilde{g}}(X)^{{2\over n}}\geq 
2^{{2\over n}}\(\min_{W} R_{\hat{g}}\)
\Vol_{\hat{g}}(W)^{{2\over n}}
\\
\\
&\geq & \displaystyle
2^{{2\over n}}(R_{\bar{g}}-\epsilon) 
(\Vol_{\bar{g}}(W)^{{2\over n}}-\epsilon) \geq 2^{{2\over n}}
R_{\bar{g}}  \Vol_{\bar{g}}(W)^{{2\over n}} - K\epsilon
\\
\\
&= & \displaystyle
2^{{2\over n}}Y_{\bar{C}}(W,M,;C) - K\epsilon.
\end{array}
$$
Here the constant $K>0$ is independent of $\epsilon$.  Thus
\vspace{2mm}

\noindent
\hspace*{30mm}$Y(X,M_1\sqcup (-M_1); C_1\sqcup C_1) \geq 2^{{2\over
n}}Y_{\bar{C}}(W,M,;C)$.
\end{Proof}
When the manifolds $\p W= M = M_{0}$, and $M_{1}$ is empty, (in this case,
the boundary of $X = W\cup_{M}(-W)$ is empty)  the following holds.
\begin{Corollary}\label{[Corollary 4.12]}$\mbox{\ }$
\begin{enumerate}
\item[{\bf (1)}]   
If  $Y_{\bar{C}}(W, M; C)\leq 0$, \ \ then \ \ 
$
    2^{{2\over n}} Y_{\bar{C}}(W, M; C)  \leq  Y(X).
$ 
\item[{\bf (2)}]   If  $Y(W, M; C)\leq 0$, \ \ \  then \ \
$    
 2^{{2\over n}} Y(W, M; C)  \leq  Y(X). 
$ 
\item[{\bf (3)}]   If  $Y(W, M)\leq 0$, \ \ \ \ \ \ \ then \ \ 
$
   2^{{2\over n}} Y(W, M) \leq Y(X). 
$
\end{enumerate}
\end{Corollary}
Now let $N$ be a smooth compact closed manifold, $\dim N = n$, and
$D^n\subset N$ be an embedded disk. Let $W= N\setminus \iint (D^n)$, with
$\p W= S^{n-1}$. 
\vspace{2mm}

\noindent
{\bf Remark.} Let $N$ be an enlargeable manifold (see \cite{GL2}).
Then the manifold $N\# (-N)= W\cup_{S^{n-1}}(-W)$ is also
enlargeable. Thus $Y(N\# (-N))\leq 0$.
\begin{Corollary}\label{enlarg}
Let $N$ be an enlargeable compact closed manifold. Then $$2^{{2\over
n}} Y(N\setminus \iint (D^n), S^{n-1})\leq Y(N\# (-N)).$$ In particular,
$Y(N\setminus \iint (D^n), S^{n-1})\leq 0$.
\end{Corollary}
{\bf Examples.} Let $T^n$ be a torus, ${\mathbf H}^n$ be a hyperbolic
space, and $\Gamma$ be a discrete group acting freely on ${\mathbf H}^n$,
so that ${\mathbf H}^n/\Gamma$ is a compact manifold. Then we have
$$
Y(T^n\setminus \iint (D^n), S^{n-1})\leq 0, \ \ \ 
Y(({\mathbf H}^n/\Gamma)\setminus \iint (D^n), S^{n-1})\leq 0. 
$$
{\bf Remark.} Let  $W_{j}$  be a compact smooth $n$-manifold 
with boundary  $\p W_j= M_{j}$,  for  $j = 1, 2$.  
Let  $(W,M) = (W_{1}, M_{1}) \sqcup (W_{2}, M_{2})$  
be the disjoint union of  $W_{1}$  and  $W_{2}$.  
Let  $C = C_{1} \sqcup C_{2}$  be a conformal class 
on  $M_{1} \sqcup M_{2}$.  
Similar to the case of closed manifolds, 
we can show that the same inequalities 
as those of \cite[Corollary 1.11]{Kobayashi1}  hold    
for the relative Yamabe invariants  
$Y(W, M; C)$  and  $Y(W_{j}, M_{j}; C_{j})$  ($j = 1, 2$).  
\section{Notes on moduli spaces}\label{open}
{\bf \ref{open}.1. Moduli space of positive scalar curvature metrics.}
Let $M$ be a closed manifold admitting a positive scalar curvature
metric. Then one has the space of psc-metrics
$$
\Riem^+(M) = \{ g\in \Riem(M) \ | \ R_g >0\}.
$$
It is known that this space has, in general, many connective
components, and that its homotopy groups are nontrivial. For
simplicity we assume that $M$ is an oriented manifold. We denote
$\Diff_+(M)$ the group of diffeomorphisms preserving the
orientation. Then the group $\Diff_+(M)$ naturally acts on the space
of metrics by pulling back a metric via a diffeomorphism. Clearly this
action preserves the space $\Riem^+(M)$. Then the moduli space of
psc-metrics is defined as ${\cal M}^+(M) = \Riem^+(M)/\Diff_+(M)$. It
is very challenging problem to describe (in some reasonable terms) the
topology of the moduli space ${\cal M}^+(M)$.  We suggest here to give
an alternative model of the moduli space of psc-metrics. First, we
suggest to start with the space ${\cal C}^+(M)$ of positive conformal
classes. Clearly there is a canonical projection map $p : \Riem^+(M)
\lra {\cal C}^+(M)$, which sends a metric $g$ to its conformal class
$[g]$. We prove the following fact.
\begin{Theorem}\label{mod1}
Let $M$ be a closed compact manifold with $\dim M \geq 3$. Then the
natural projection map $p: \Riem^+(M)\lra {\cal C}^+(M)$ is weak
homotopy equivalence.
\end{Theorem}
\begin{Proof}
We start with the following easy observation.
\begin{Lemma}\label{mod-1}
Let $C\in {\cal C}^+(M)$, and $g_0, g_1\in C$ be psc-metrics. Then
$g_0$ and $g_1$ are psc-homotopic, i.e. there exists a smooth family
$\{g(t)\}_{t\in [0,1]} \in C$ of psc-metrics with $g(0)=g_0$,
$g(1)=g_1$.
\end{Lemma}
\begin{Proof} Indeed, we have that $R_{g_0}>0$, and $g_1=u^{{4\over
n-2}}g_0$ for $u\in C^{\infty}_+(M)$ with the scalar curvature
$$
R_{g_1} = u^{-{n+2\over n-2}}\( -{4(n-1)\over (n-2)}\Delta u +
R_{g_0}u\) > 0. 
$$
Then the curve of metrics $g(t) = u(t)^{{4\over n-2}}g_0 \in C$ with
$u(t) = ut + (1-t) > 0$ is such that
$$
\begin{array}{lcl}
R_{g(t)} & = &\displaystyle
u(t)^{-{n+2\over n-2}} \( -{4(n-1)\over n-2}\Delta u(t) + R_{g_0} u(t) \)
\\
& = &\displaystyle
u(t)^{-{n+2\over n-2}} \( 
t\[-{4(n-1)\over n-2}\Delta u +  R_{g_0} u\] + (1-t)R_{g_0}  \)
\\
& = &\displaystyle
u(t)^{-{n+2\over n-2}}
\(
t R_{g_1} u^{{n+2\over n-2}} + (1-t) R_{g_0} 
\) > 0
\end{array}
$$ 
since the functions $R_{g_1} u^{{n+2\over n-2}}$ and $R_{g_0}$ are
both positive.  
\end{Proof}
Now let $P(C)= \{ g\in C \ | \ R_g >0 \ \}$. Clearly $C\cong
C^{\infty}_+(M)$ is a convex set.
\begin{Lemma}\label{mod-2}
The subset $P(C)\subset C$ is a convex and contractible set.
\end{Lemma}
\begin{Proof} First we check that $P(C)$ is convex. Indeed, let
$\check{g}\in C$ be a Yamabe metric, with $R_{\check{g}}={\mathrm
c}{\mathrm o}{\mathrm n}{\mathrm s}{\mathrm t}.  >0$. Then for any
$g\in P(C)$ there exists a unique function $u\in C^{\infty}_+(M)$ so
that $g = u^{{4\over n-2}}\check{g}$. Thus we identify $P(C)$ with the
following subspace of positive smooth functions
$$
P(C)\cong\{ u\in C^{\infty}_+(M) \ | \ -\Delta u + {(n-2)\over 4(n-1)
}R_{\check{g}}u > 0 \} .
$$
A homotopy $F_t : P(C) \longrightarrow P(C)$ given by
$$
F_t(g) = u^{{4\over n-2}}(t)\check{g} \ \ \mbox{with} \ \ 
u(t) = ut + (1-t)
$$
is well defined by Lemma \ref{mod-1}. Futhermore, $F_1 = Id$, and $F_0
$ sends the set $P(C)$ to a single point $\check{g}\in P(C)$.
Therefore $P(C)$ is convex, and since $P(C)$ is a subspace of the
convex space $C^{\infty}_+(M)$, it is and contractible.
\end{Proof}
We notice that both spaces $\Riem^+(M)$ and ${\cal C}^+(M)$ have
homotopy types of $CW$-complexes. Thus we can assume (up to homotopy
equivalence) that $p: \Riem^+(M) \longrightarrow {\cal C}^+(M)$ is a
fibration. Since $p^{-1}(C)$ is contractible for any conformal class
$C$, we obtain that $p$ induces isomorphism in homotopy groups $p_* :
\pi_k(\Riem^+(M))\cong \pi_k({\cal C}^+(M))$.
\end{Proof}
Thus in the homotopy category one does not loose any information by
replacing the space $\Riem^+(M)$ by the space of positive conformal
classes ${\cal C}^+(M)$. 
\vspace{2mm}

\noindent
The space ${\cal C}(M)$ is the orbit space of the action (left
multiplication) of the group $C_+^{\infty}(M)$ on the space of metrics
$\Riem(M)$. It is convenient to refine this construction (as it is
done in \cite{Morava}) for manifolds with a base point.
\vspace{2mm}

\noindent
Let $x_0\in M$ be a base point. We consider the following subspace of
$C_+^{\infty}(M)$:
$$
C_{+,x_0}^{\infty}(M) = \{ u\in C_+^{\infty}(M) \ | \ u(x_0)= 1 \ \}.
$$ 
Then let ${\cal C}_{x_0}(M)$ be the orbit space of the induced action
of $C_{+,x_0}^{\infty}(M)$ on $\Riem(M)$. Clearly there is a canonical
map $p_1: {\cal C}_{x_0}(M) \lra {\cal C}(M)$ which is a homotopy
equivalence since $p_1^{-1}(C) \cong \R$. 
Let
$$
{\cal C}_{x_0}^+(M) = p_1^{-1} \( {\cal C}^+(M)\).
$$
To construct an appropriate moduli space we assume that $M$ is a
connected manifold, and consider the following subgroup of the
diffeomorphism group $\Diff_+(M)$:
$$
\Diff_{x_0,+}(M)=\{\phi\in \Diff_+(M) \ | \ \phi(x_0)=x_0, \ \ d\phi=Id :
TM_{x_0} \rightarrow TM_{x_0} \}
$$
The group $\Diff_{x_0,+}(M)$ inherits the action on the spaces ${\cal
C}(M)$ and ${\cal C}_{x_0}(M) $. It is easy to prove the group
$\Diff_{x_0,+}(M)$ acts freely on the space ${\cal C}_{x_0}(M)$
(perhaps, it is important that $M$ is connected).
\vspace{2mm}

\noindent
Clearly the space ${\cal C}_{x_0}^+(M) $ of positive conformal classes
is invariant under this action. We define the moduli space ${\cal
M}^{+}_{x_0,\conf}(M)$ of positive conformal structures as the orbit space
of the action of $\Diff_{x_0,+}(M)$ on ${\cal C}_{x_0}^+(M)$.  One
obtains the diagram of Serre fiber bundles 
$$
\begin{diagram}
\node{{\cal C}_{x_0}^+(M)} 
\arrow{s,l}{\pi^+}
\arrow[2]{e,t}{i}
\node[2]{{\cal C}_{x_0}(M)}
\arrow{s,l}{\pi}
\\
\node{{\cal M}^{+}_{x_0,\conf}(M)} 
\arrow[2]{e,t}{i^{+}}
\node[2]{B\Diff_{x_0,+}(M)}
\end{diagram}
$$
Here $B\Diff_{x_0,+}(M)$ is the classifying space of the group
$\Diff_{x_0,+}(M)$ which we identify with the orbit space ${\cal
C}_{x_0}(M)/ \Diff_{x_0,+}(M)$ (since the action is free, and the space
${\cal C}_{x_0}(M)$ is contractible). We address the following problem.
\begin{Problem} What is the rational homotopy type of the space
${\cal M}^{+}_{x_0,\conf}(M)$?
\end{Problem}
{\bf \ref{open}.2. Conformal isotopy and concordance.} It is
well-known that isotopic psc-metrics are concordant, see \cite{GL1} and
\cite{Gajer}. It is still not known if the converse is true; (we quote
\cite{RS}) ``indeed, there is no known method to distinguish between
isotopy classes of positive scalar curvature which is not based on
distinguishing concordance classes.'' We would like to address the
``conformal analogue'' of this problem.
\vspace{2mm}

\noindent
Let $C_0, C_1\in {\cal C}^+(M)$ be two positive conformal
classes.  One defines an isotopy of positive conformal classes in the
obvios way.  We say that the conformal classes $C_0$ and $C_1$
are {\it conformally concordant} if
$$
Y(M\times [0,1], M\times \{0,1\}; C_0\sqcup
C_1)>0.
$$
Theorem \ref{ThI} implies the following result:
\begin{Corollary}\label{CorI}
Conformal concordance is an equivalence relation on ${\cal C}^+(M)$.
\end{Corollary}
We would like to spell out the following conjecture:
\begin{Conjecture}\label{ABConj2}
Let $M$ be a closed compact manifold admitting a psc-metric, $n\geq
5$. If $C_0,C_1\in {\cal C}^+(M)$ are conformally concordant, the
the classes $C_0,C_1$ are isotopic in ${\cal C}^+(M)$.
\end{Conjecture}
\vspace{2mm}

\noindent
{\bf \ref{open}.3. Conformal cobordism.} Once we would like to
describe the whole world of manifolds equipped with psc-metrics, we
are led to a concept of cobordism.  Two manifolds $(M_0,g_0)$,
$(M_1,g_1)$ with psc-metrics $g_0, g_1$, are said to be psc-cobordant
if there exists a manifold $(W,\bar{g})$ with $\p W= M_0\sqcup
(-M_1)$, and a psc-metric $\bar{g}$, so that:
\begin{enumerate}
\item[{\bf (1)}] $\bar{g}|_{M_j}= g_j$, $j=1,2$,
\item[{\bf (2)}] $\bar{g} = g_j + dr^2$ near $M_j$.
\end{enumerate}
We emphasize that the metric $\bar{g}$ must be a product metric near
the boundary. The psc-cobordims was used in several papers \cite{BG},
\cite{Gajer}, \cite{Hajduk}, \cite{Stolz1}. For instance, S. Stolz
described an adequate psc-cobordism category where given manifold $M$
fits in (see \cite{Stolz1}). This category is determined by the
fundametal group $\pi_1(M)$ and the first two Stiefel-Whitney classes
of $M$. 
\vspace{2mm}

\noindent
We define the conformal analogue of the psc-cobordism relation by
means of the relative Yamabe invariant. Let $(M_0,C_0)$, $(M_1,C_1)$
be two manifolds equipped with positive conformal classes. We call
such manifolds {\sl positive conformal manifolds}.  Then $(M_0,C_0)$,
$(M_1,C_1)$ are {\sl conformally cobordant} if there is a manifold
$W$, with $\p W= M_0\sqcup (-M_1)$, and such that the relative Yamabe
invariant
$$
Y(W,M_0\sqcup (-M_1);C_0\sqcup C_1)>0.
$$
Again, Theorem \ref{ThI} implies the following
result:
\begin{Corollary}\label{CorII} Conformal cobordism is
an equivalence relation on the category of positive conformal
manifolds.
\end{Corollary}
{\bf Remark.} The definition of the conformal cobordism may be
essentially refined in the way suggested by S. Stolz \cite{Stolz1}.
This leads to the corresponding conformal cobordism groups. We are
studying these  cobordism groups in another paper.

\newpage
 
%
%
%


\begin{thebibliography}{30}
\begin{small}
\bibitem{Besse} A. Besse, Einstein manifolds, Springer, Berlin, 1987.
\bibitem{BG} B. Botvinnik, P. Gilkey, The eta invariant and metrics of
	positive scalar curvature,  Math. Ann. 302 (1995), 507--517.
\bibitem{Chavel} I. Chavel, Eigenvalues in Riemannian geometry, Pure
	and Applied Mathematics, 115. Academic Press,  1984.
\bibitem{Cherrier} P. Cherrier, Probl\'{e}mes de Neumann non
	lin\`{e}aires sur les vari\`{e}t\`{e}s riemanniennes,
	J. Funct. Anal. 57 (1984),  154--206.
\bibitem{Escobar1} J. Escobar, The Yamabe problem on manifolds with
        boundary, J. Diff. Geom. 35 (1992), 21-84.
\bibitem{Gajer} P. Gajer, Concordances of metrics of positive scalar
	curvature, Pacific J.  Math. 157 (1993),  257--268
\bibitem{GL1} M. Gromov, H. B. Lawson, The classification
        of simply connected manifolds of positive scalar curvature,
        Ann. of Math. 111 (1980), 423-434.
\bibitem{GL2} M.  Gromov, H. B. Lawson, Positive scalar curvature and
	the Dirac operator on complete Riemannian manifolds,
	Publ. Math. Inst. Hautes \'{E}tudes Sci.  No. 58, (1983),
	83--196 (1984).
\bibitem{Hajduk} B. Hajduk, On the obstruction group to existence of
	Riemannian metrics of positive scalar curvature, Global
	differential geometry and global analysis (Berlin, 1990),
	62--72, Lecture Notes in Math. 1481, Springer, Berlin, 1991.
\bibitem{Kobayashi1} O. Kobayashi, Scalar curvature of a metric with
	unit volume, Math. Ann. 279 (1987), 253-265.
\bibitem{Morava} J. Morava, H. Tamanoi, A vanishing theorem for the
	conformal anomaly in dimension
	$>2$, Proc. Amer. Math. Soc. 100 (1987), 767--774.
\bibitem{RS} J. Rosenberg, S. Stolz, Metric of positive scalar
	curvature and connection with surgery, to appear.
\bibitem{Stolz1} S. Stolz, Concordance classes of positive scalar
	curvature metrics, to appear.
\end{small} 
\end{thebibliography}
\end{document}